\documentstyle[11pt]{article}
  \newlength{\standardunitlength}
\newtheorem{cor}{Corollary} \newtheorem{lemma}{Lemma}
\newtheorem{theorem}{Theorem} \newtheorem{prop}{Proposition}
\newenvironment{proof}{\noindent {\sc Proof:}}{$\Box$ \vspace{2 ex}}

\begin{document}

\begin{center}
Cycle Indices for the Finite Classical Groups
\end{center}

\begin{center}
By Jason Fulman
\end{center}

\begin{center}
Dartmouth College
\end{center}

\begin{center}
Jason.E.Fulman@Dartmouth.Edu
\end{center}

\begin{abstract}
	This paper defines and develops cycle indices for the finite classical
groups. These tools are then applied to study properties of a random matrix chosen
uniformly from one of these groups. Properties studied by this technique will
include semisimplicity, regularity, regular semisimplicity, the characteristic
polynomial, number of Jordan blocks, and average order of a matrix.
\end{abstract}

\section{Introduction and Background} \label{Intro}

	Polya \cite{Polya}, in a landmark paper on combinatorics,
introduced the cycle index of the symmetric groups. This can be written as
follows. Let $a_i(\pi)$ be the number of $i$-cycles of $\pi$. The Taylor
expansion of $e^z$ and the fact that there are $\frac{n!} {\prod_i a_i! i^{a_i}}$
elements of $S_n$ with $a_i$ $i$-cycles, yield the following theorem.

\begin{theorem} (Polya \cite{Polya})

\[ \sum_{n=0}^{\infty} \frac{u^n}{n!} \sum_{\pi \in S_n} \prod_i
x_i^{a_i(\pi)} = \prod_{m=1}^{\infty} e^{\frac{x_mu^m}{m}}. \]

\end{theorem}

	The Polya cycle index has been a key tool in understanding what a
typical permutation $\pi \in S_n$ "looks like". It is useful for studying
properties of a permutation which depend only on its cycle structure. Here
are some examples of theorems which can be proved using the cycle
index. Lloyd and Shepp \cite{She} showed that for any $i< \infty$, the
joint distribution of $(a_1(\pi),\cdots,a_i(\pi))$ for $\pi$ chosen
uniformly in $S_n$ converges to independent (Poisson($1$), $\cdots$,
Poisson($\frac{1}{i}$)) as $n \rightarrow \infty$. Goncharov $\cite{Gon}$
proved that the number of cycles in a random permutation is asymptotically
normal with mean $log\ n$ and standard deviation $(log \
n)^{\frac{1}{2}}$. Goh and Schmutz $\cite{Go1}$ proved that
if $\mu_n$ is the average order of an element of $S_n$, then:

\[ log \ \mu_n = C \sqrt{\frac{n}{log \ n}} + O(\frac{\sqrt{n} log log
n}{log n}) \]

	where $C = 2.99047...$. Letting $L_1$ denote the
length of the longest cycle of a permutation, Goncharov $\cite{Gon}$ showed that
$\frac{L_1}{n}$ has a limit distribution. Lloyd and Shepp $\cite{She}$ computed
the moments of this limit distribution. For instance, letting $E_n$ signify
expectation in $S_n$, they proved that:

	\[ \lim_{n \rightarrow \infty} E_n (\frac{L_1}{n}) =
\int_0^{\infty} exp [-x - \int_x^{\infty} \frac{e^{-y}}{y} dy] dx =
0.62432997... \]

	For more on the Polya cycle index, including a wonderful probabilistic
interpretation due to Lloyd and Shepp, see Chapter 1 of Fulman \cite{fulthesis}.

	Kung \cite{Kung} and Stong \cite{St1} developed a cycle index for
the finite general linear groups, with an eye toward studying properties of a
random element of $GL(n,q)$. Chapter 1 of Fulman \cite{fulthesis} gives motivation
(for instance from the theory of random number generators) for why one might want
to study random matrices. Let us review the Kung/Stong cycle index, as one of the
goals of this paper is to carry out a similar undertaking for the other finite
classical groups. First it is necessary to understand the conjugacy classes of
$GL(n,q)$. As is explained in Chapter 6 of Herstein $\cite{Her}$, an element
$\alpha \in GL(n,q)$ has its conjugacy class determined by its rational canonical
form (this is slightly different from Jordan canonical form, which only works
over algebraically closed fields). This form corresponds to the following
combinatorial data. To each monic non-constant irreducible polynomial $\phi$ over
$F_q$, associate a partition (perhaps the trivial partition) $\lambda_{\phi}$ of
some non-negative integer $|\lambda_{\phi}|$ into parts $\lambda_1 \geq \lambda_2 \geq \cdots$. Let $m_{\phi}$ denote the
degree of $\phi$. The only restrictions necessary for this data to
represent a conjugacy class are: conjugacy class are:

\begin{enumerate}

\item $|\lambda_z| = 0$
\item $\sum_{\phi} |\lambda_{\phi}| m_{\phi} = n$

\end{enumerate}

	An explicit representative of this conjugacy class may be given as
follows. Define the companion matrix $C(\phi)$ of a polynomial
$\phi(z)=z^{m_{\phi}} + \alpha_{m_{\phi}-1} z^{m_{\phi}-1} + \cdots +
\alpha_1 z + \alpha_0$ to be:

\[ \left( \begin{array}{c c c c c}
		0 & 1 & 0 & \cdots & 0 \\
		0 & 0 & 1 & \cdots & 0 \\
		\cdots & \cdots & \cdots & \cdots & \cdots \\
		0 & 0 & 0 & \cdots & 1 \\
		-\alpha_0 & -\alpha_1 & \cdots & \cdots & -\alpha_{m_{\phi}-1}
	  \end{array} \right) \]
	
	Let $\phi_1,\cdots,\phi_k$ be the polynomials such that
$|\lambda_{\phi_i}|>0$. Denote the parts of $\lambda_{\phi_i}$ by
$\lambda_{\phi_i,1} \geq \lambda_{\phi_i,2} \geq \cdots $. Then a matrix
corresponding to the above conjugacy class data is:

\[ \left( \begin{array}{c c c c}
		R_1 & 0 & 0 &0 \\
		0 & R_2 & 0 & 0\\
		\cdots & \cdots & \cdots & \cdots \\
		0 & 0 & 0 & R_k
	  \end{array} \right) \]

	where $R_i$ is the matrix:

\[ \left( \begin{array}{c c c}
		C(\phi_i^{\lambda_{\phi_i,1}}) & 0 & 0  \\
		0 & C(\phi_i^{\lambda_{\phi_i,2}}) & 0 \\
		0 & 0 & \cdots 
	  \end{array} \right) \]

	For example, the identity matrix has $\lambda_{z-1}$ equal to
$(1^n)$ and an elementary reflection with $a \neq 0$ in the $(1,2)$
position, ones on the diagonal and zeros elsewhere has $\lambda_{z-1}$
equal to $(2,1^{n-2})$.

	As was true for the symmetric groups, many algebraic
properties of a matrix $\alpha$ can be stated in terms of the data
parameterizing its conjugacy class. For example, the characteristic polynomial of $\alpha \in GL(n,q)$ is equal
to $\prod_{\phi} \phi^{|\lambda_{\phi}(\alpha)|}$. Section \ref{REG} will give further examples.

	Kung \cite{Kung} and Stong \cite{St1} developed a cycle index for the finite
general linear groups, analogous to Polya's cycle index for the symmetric groups.
It can be described as follows. Let $x_{\phi,\lambda}$ be variables corresponding
to pairs of polynomials and partitions. Then the cycle index for $GL(n,q)$ is
defined by

	\[ Z_{GL(n,q)} = \frac{1}{|GL(n,q)|} \sum_{\alpha \in GL(n,q)}
\prod_{\phi \neq z} x_{\phi,\lambda_{\phi}(\alpha)} \]
	
	It is also helpful to define a quantity
$c_{GL,\phi,q^{m_{\phi}}}(\lambda)$. If $\lambda$ is the empty partition, set
$c_{GL,\phi,q^{m_{\phi}}}(\lambda)=1$. Using the standard notation for
partitions, let $\lambda$ have $m_i$ parts of size $i$. Write:

\[ d_i= m_1 1 + m_2 2 + \cdots + m_{i-1}(i-1) + (m_i + m_{i+1} +
\cdots + m_j) i \]

	Then define:

\[ c_{GL,\phi,q^{m_{\phi}}}(\lambda) = \prod_i \prod_{k=1}^{m_i} (q^{m_{\phi}d_i} -
q^{m_{\phi}(d_i-k)}) \]

	Following Kung $\cite{Kung}$, Stong $\cite{St1}$ proves the
factorization:

	\[ 1+\sum_{n=1}^{\infty} Z_{GL(n,q)} u^n = \prod_{\phi \neq z}
\sum_{\lambda} x_{\phi,\lambda}
\frac{u^{|\lambda|m_{\phi}}}{c_{GL,\phi,q^{m_{\phi}}}(\lambda)} \]

	The terms $c_{GL,\phi,q^{m_{\phi}}}$ in the Kung-Stong
cycle index appear difficult to work with, and Section $\ref{GLCYC}$ will give several useful rewritings. Nevertheless, Stong $\cite{St1}, \cite{St2}$ and Goh and
Schmutz $\cite{Go2}$ have used this cycle index successfully to study properties of random elements of
$GL(n,q)$. For instance,  Stong \cite{St1} 
showed that the number of Jordan blocks of a random element of $GL(n,q)$
has mean and variance $log(n)+O(1)$, and Goh and Schmutz \cite{Go2} proved convergence
to the normal distribution. Stong $\cite{St1}$ also obtained asymptotic (in $n$ and $q$) estimates for the 
chance that an element of $GL(n,q)$ is a vector space
derangement (i.e. fixes only the origin),  the number of elements of $GL(n,q)$ satisfying a fixed
polynomial equation, and the chance that all polynomials appearing in the rational
canonical form of $\alpha \in GL(n,q)$ are linear. Stong \cite{St2} studied the average order of a matrix.

	The structure of this paper is as follows. Section \ref{GLCYC} gives useful
rewritings of the cycle index for the general linear groups. Section \ref{REG}
applies the Kung/Stong cycle index to obtain exact formulas for
the $n \rightarrow \infty$ limit of the chance that an element of
$GL(n,q)$ or $Mat(n,q)$ is semisimple, regular, or regular semisimple. Section
\ref{OTHER} uses work of Wall on the conjugacy classes of the unitary,
symplectic, and orthogonal groups to obtain cycle indices for these groups.
Section \ref{BIGQ} uses the theory of algebraic groups to study the $q
\rightarrow \infty$ limit of these cycle indices. Finally, Section
\ref{APPLICATIONS} applies the cycle indices of Section \ref{OTHER} to study the
characteristic polynomial, number of Jordan blocks, and average order of an
element of a finite classical group. There is much more applied work to be done
using these cycle indices; the results of Sections \ref{REG} and
\ref{APPLICATIONS} are meant only as a start. Section \ref{SUGGESTIONS} gives
some suggestions for future research.

	The work in this paper is taken from the author's Ph.D. thesis \cite{fulthesis}, done under the supervision of Persi Diaconis. There are three papers which should be considered companions to this one. Fulman \cite{fulmacdonald} connects the cycle indices with symmetric function theory and
exploits this connection to develop probabilistic algorithms with group theoretic
meaning. Fulman \cite{fulprob} applies the algorithms of Fulman
\cite{fulmacdonald} to prove group theoretic results and can be read with no
knowlegde of symmetric functions. Finally, Fulman \cite{fulRogers} establishes an
appearance of the Rogers-Ramanujan identities in the finite general linear
groups; this is relevant to Section \ref{REG} of this paper.

\section{Cycle Index of the General Linear Group: Useful Rewritings}
\label{GLCYC} 

	Although the ultimate purpose of this section is to rewrite the cycle index
of the general linear groups in useful ways, we begin with some remarks about the
cycle index which were excluded from the brief discussion in the introductory
section. First recall from the introduction that:

\[ 1+\sum_{n=1}^{\infty} \frac{u^n}{|GL(n,q)|} \sum_{\alpha \in GL(n,q)}
\prod_{\phi \neq z} x_{\phi,\lambda_{\phi}(\alpha)}  = \prod_{\phi \neq z}
\sum_{\lambda} x_{\phi,\lambda}
\frac{u^{|\lambda|m_{\phi}}}{c_{GL,\phi,q^{m_{\phi}}}(\lambda)} \]

	where

\[ c_{GL,\phi,q^{m_{\phi}}}(\lambda) = \prod_i \prod_{k=1}^{m_i} (q^{m_{\phi}d_i} -
q^{m_{\phi}(d_i-k)}) \]

	and

\[ d_i= m_1 1 + m_2 2 + \cdots + m_{i-1}(i-1) + (m_i + m_{i+1} +
\cdots + m_j) i .\]

{\bf Remarks}

\begin{enumerate}

\item The fact that the cycle index factors comes from the fact in Kung \cite{Kung} that if $\alpha$ has data $\lambda_{\phi}(\alpha)$, then
the conjugacy class of $\alpha$ in $GL(n,q)$ has size:

\[ \frac{|GL(n,q)|}{\prod_{\phi}
c_{GL,\phi,q(\lambda_{\phi}(\alpha))}} \]

	The following example should make this formula seem more real. A
transvection in $GL(n,q)$ is defined as a determinant 1 linear map whose
pointwise fixed space is a hyperplane.  For instance the matrix with $a
\neq 0$ in the $(1,2)$ position, ones on the diagonal and zeros elsewhere
is a transvection. The transvections generate $SL(n,q)$ (e.g. Suzuki
$\cite{Suzuki}$) and are useful in proving the simplicity of the projective
special linear groups.

	Transvections can be counted directly. Let $V$ be an $n$ dimensional vector space
over
$F_q$ with dual space $V^*$. It is not hard to see that the action of any
transvection $\tau$ is of the form:

\[ \tau(\vec{x}) = \vec{x} + \vec{a} \psi(\vec{x}) \]

	where $\vec{a} \in V$ is a non-0 vector and $\psi \in V^*$ is
a non-0 linear form on $V$ which vanishes on $\vec{a}$.  It readily follows that the number
of transvections is:

\[ \frac{(q^n-1)(q^{n-1}-1)}{q-1} \]

	Transvections can also be counted using Kung's class size formula. It is not hard to show that an element $\alpha \in GL(n,q)$ is a transvection if
and only if $\lambda_{z-1}(\alpha)=(2,1^{n-2})$ and
$|\lambda_{\phi}(\alpha)|=0$ for all $\phi \neq z-1$. This follows from the
fact that a transvection has all eigenvalues 1 and from a lemma in Fulman
\cite{fulthesis}, which says that the dimension of the fixed space of
$\alpha$ is the number of parts of the partition
$\lambda_{z-1}(\alpha)$. Thus all transvections in $GL(n,q)$ are
conjugate and Kung's formula shows that the size of the conjugacy
class in $GL(n,q)$ corresponding to $\lambda_{z-1}=(2,1^{n-2})$ is:

\[ \frac{|GL(n,q)|}{c_{GL,z-1,q}(2,1^{n-2})} = \frac{(q^n-1)
(q^{n-1}-1)}{q-1}. \]

\item The orbits of $GL(n,q)$ on $Mat(n,q)$ (all $n*n$ matrices)
under conjugation are again parameterized by the data
$\lambda_{\phi}$. However the parition $\lambda_{\phi}$ need not have size zero. This leads to the factorization:

\[ 1+\sum_{n=1}^{\infty} \frac{u^n}{|Mat(n,q)|} \sum_{\alpha \in Mat(n,q)}
\prod_{\phi \neq z} x_{\phi,\lambda_{\phi}(\alpha)}  = \prod_{\phi}
\sum_{\lambda} x_{\phi,\lambda}
\frac{u^{|\lambda|m_{\phi}}}{c_{GL,\phi,q^{m_{\phi}}}(\lambda)} \]

	This factorization will be used in Section \ref{REG}, and shows that any successful application of cycle index theory to $GL(n,q)$ can be carried over to $Mat(n,q)$, and conversely.

\item Many of Stong's applications of the general linear group cycle index made use of the
fact that it is possible to count the number of monic, degree $m$, irreducible
polynomials $\phi \neq z$ over a finite field $F_q$. Letting $I_{m,q}$ denote the
number of such polynomials, the following lemma is well-known. Let $\mu$ be the
usual Moebius function of elementary number theory.

\begin{lemma} \label{irred}
\[ I_{m,q} = \frac{1}{m} \sum_{k|m} \mu(k) (q^{\frac{m}{k}}-1) \]
\end{lemma}

\begin{proof}
	From Hardy and Wright $\cite{HW}$, the number of monic, degree
$m$, irreducible polynomials with coefficients in $F_q$ is:

	\[ \frac{1}{m} \sum_{k|m} \mu(k) q^{\frac{m}{k}} \]

	Now use the fact that $\sum_{k|m} \mu(k)$ is $1$ if $m=1$ and $0$
otherwise.
\end{proof}

	Similar counting results, for special types of polynomials relevant to the other finite classical groups, will be obtained in Section \ref{OTHER}.

\end{enumerate}

	Now, we give two rewritings of the seemingly horrible quantitiy $c_{GL,\phi,q^{m_{\phi}}}$ which appeared in the Kung/Stong cycle index for the general linear groups. Recall that if $\lambda$ is a partition of $n$ into parts $\lambda_1 \geq \lambda_2 \cdots$, that
$m_i(\lambda)$ is the number of parts in the partition equal to $i$. Let
$\lambda_i'= m_i(\lambda) + m_{i+1}(\lambda) + \cdots$ be the $i$th part of the
partition dual to $\lambda$. Also, $(\frac{1}{q})_r$ will denote $(1-\frac{1}{q})
\cdots (1-\frac{1}{q^r})$.

\begin{theorem} \label{Rewrite}

\begin{eqnarray*}
c_{GL,\phi,q^{m_{\phi}}}(\lambda) & = & q^{2m_{\phi} [\sum_{h<i} h
m_h(\lambda) m_i(\lambda) + \frac{1}{2} \sum_i (i-1) m_i(\lambda)^2]}
\prod_i |GL(m_i(\lambda),q^{m_{\phi}})|\\
& = & q^{m_{\phi} [\sum_i
(\lambda'_i)^2]} \prod_{i} (\frac{1}{q^{m_{\phi}}})_{m_i(\lambda)}
\end{eqnarray*}

\end{theorem}

\begin{proof}
	For both equalities assume that $m_{\phi}=1$, since the
result will be proved for all $q$ and one could then substitute
$q^{m_{\phi}}$ for $q$.

	For the first equality, it's easy to see that the factors of
the form $q^r-1$ are the same on both sides, so it suffices to look at
the powers of $q$ on both sides. The power of $q$ on the left-hand
side is $\sum_i [d_im_i(\lambda) - {m_i(\lambda) \choose 2}]$ and the
power of $q$ on the right-hand side is $\sum_i [i m_i(\lambda)^2
-{m_i(\lambda) \choose 2} + \sum_{h<i} hm_h(\lambda)
m_i(\lambda)]$. Thus it is enough to show that:

\[ \sum_i d_i = \sum_i [im_i(\lambda) + 2 \sum_{h<i} hm_h(\lambda)] \]

	This equality follows from the observation that:

\[ d_i = [\sum_{h<i} hm_h(\lambda)] + im_i(\lambda) + [\sum_{i<k}
im_k(\lambda)] \]

	For the second equality of the theorem, write
$(\frac{1}{q})_{m_i(\lambda)}$ as $\frac{|GL(m_i(\lambda),q)|}
{q^{m_i(\lambda)^2}}$. Comparing powers of $q$ reduces us to proving
that:

\[ \sum_i (\lambda_i')^2 = \sum_i [im_i(\lambda) + 2 \sum_{h<i}
hm_h(\lambda)] m_i(\lambda) \]

	This last equation follows quickly after substituting
$\lambda_i'=m_i(\lambda)+m_{i+1}(\lambda)+\cdots$.
\end{proof} 

	In fact, there is a third (and extremely useful) rewriting of
$c_{GL,\phi,q^{m_{\phi}}}(\lambda)$ in terms of $P_{\lambda}(x_i;t)$, the
Hall-Littlewood symmetric functions. Let $n(\lambda) = \sum_i (i-1) \lambda_i =
\sum_i {\lambda_i'
\choose 2}$. Theorem 3 of Fulman
\cite{fulthesis} shows that:

\[  c_{GL,\phi,q^{m_{\phi}}}(\lambda)   = 
\frac{q^{m_{\phi}
n(\lambda)}}{P_{\lambda}(\frac{1}{q^{m_{\phi}}},\frac{1}{q^{2m_{\phi}}},\cdots;\frac{1}{q^{m_{\phi}}})} \]

	The proof is omitted in the present treatment, as this article will not make use of the connections with symmetric function theory developed in Fulman
\cite{fulthesis}.

	Lemma \ref{Sto} is due to Stong \cite{St1} and will be useful in future
sections. It was proved using an identity of Goldman and Rota and the fact
that the number of unipotent (all eigenvalues equal to 1) elements of $GL(n,q)$ is
$q^{n(n-1)}$. Fulman
\cite{fulprob} gives a probabilstic proof.

\begin{lemma} \label{Sto}

\[ \sum_{\lambda} \frac{1}{c_{GL,\phi,q^{m_{\phi}}}(\lambda)} =
\prod_{r=1}^{\infty} (\frac{1}{1-\frac{u^{m_{\phi}}}{q^{rm_{\phi}}}}) \]

\end{lemma} 

\section{Application: Counting Semisimple, Regular, and Regular Semisimple
Matrices} \label{REG}

	This section will use the cycle index of the general linear groups to obtain
asymptotic formulas for the chance that an element of $GL(n,q)$ or $Mat(n,q)$ is
semisimple, regular, or regular semisimple. The strategy for applying the cycle
index to studying a property of a matrix which depends only on its conjugacy
class will always be the same: first translate the property into a statement
about the partitions $\lambda_{\phi}$ parameterizing the conjugacy class, then
manipulate the cycle index accordingly.

	The following three lemmas will prove useful for working with the cycle index.
Lemma \ref{bign} allows one to gain insight into $n
\rightarrow \infty$ asymptotics. (An alternate approach uses the
method of moments). We use the notation that $[u^n] g(u)$ is the
coefficient of $u^n$ in the Taylor expansion of $g(u)$ around $0$.

\begin{lemma} \label{bign} If $f(1)<\infty$ and $f$ has a Taylor
series around 0, then:

\[ lim_{n \rightarrow \infty} [u^n] \frac{f(u)}{1-u} = f(1) \]

\end{lemma}

\begin{proof}
	Write the Taylor expansion $f(u) = \sum_{n=0}^{\infty} a_n
u^n$. Then observe that $[u^n] \frac{f(u)}{1-u} = \sum_{i=0}^n a_i$.
\end{proof}

	Lemmas $\ref{allpoly}$ and $\ref{product}$ will be useful for manipulating
products which appear in the cycle index. Recall that $m_{\phi}$ is the degree of
a monic, polynomial
$\phi$ which is irreducible over $F_q$, the finite field of size $q$.

\begin{lemma}\label{allpoly} \[ \prod_{\phi} (1-\frac{u^{m_{\phi}}}{q^{m_{\phi}t}}) =
1-\frac{u}{q^{t-1}} \]
\end{lemma}

\begin{proof}
	Assume that $t=1$, the general case following by replacing $u$
by $\frac{u}{q^{t-1}}$. Expanding $\frac{1}{1-\frac{u^{m_{\phi}}}
{q^{m_{\phi}}}}$ as a geometric series and using unique factorization in
$F_q[x]$, one sees that the coefficient of $u^d$ in the reciprocal of the left
hand side is $\frac{1}{q^d}$ times the number of monic polynomials of degree $d$,
hence 1. Comparing with the reciprocal of the right hand side completes the proof.
\end{proof}

\begin{lemma} \label{product}

\[ \prod_{\phi \neq z} \prod_{r=1}^{\infty} (1-\frac{u^{m_{\phi}}}{q^{rm_{\phi}}})
= 1-u
\]

\end{lemma}

\begin{proof}
	Lemma \ref{allpoly} implies that:

\[ \prod_{\phi} \prod_{r \geq 1} (1-\frac{u^{m_{\phi}}}{q^{rm_{\phi}}}) =
\prod_{r \geq 1} (1-\frac{u}{q^{r-1}}) \]

	The result follows by cancelling the terms corresponding to $\phi=z$.
	Alternatively, this result can be deduced from Lemma \ref{Sto}.
\end{proof}
	
\begin{center}
{\bf Application 1: Semisimplicity}
\end{center}

	Recall that an element of $Mat(n,q)$ is said to be semisimple if it is
diagonalizable over $\bar{F_q}$, the algebraic closure of $F_q$. This application
uses the cycle index to study the $n \rightarrow \infty$ limit of the chance
that an element of $GL(n,q)$ or $Mat(n,q)$ is semisimple. The first step is to
express this condition in terms of rational canonical form.

	Recall the Jordan canonical form of a matrix (Chapter 6 of
Herstein $\cite{Her}$), which parameterizes the conjugacy classes of
$GL$ over an algebraically closed field, such as $\bar{F_q}$. This is
the same as the rational canonical form of a matrix (Section
$\ref{Intro}$), except that now the companion matrix $C(\phi_i^A)$
is conjugate to:

\[ \left( \begin{array}{c c c c}
		D(\beta) & 0 & 0 & 0 \\
		0 & D(\beta^q) & 0 &0  \\
		\cdots & \cdots & \cdots & \cdots \\
		0 & 0 & 0 & D(\beta^{q^{m_{\phi}-1}})
	  \end{array} \right) \]

	where $\beta, \cdots, \beta^{q^{m_{\phi}-1}}$ are the roots of
$\phi$ and $D(\gamma)$ is the $A*A$ matrix:

\[ \left( \begin{array}{c c c c c}
		\gamma & 1  & 0 & 0 & 0\\
		 0 & \gamma & 1 & 0 & 0 \\
		 \cdots & \cdots & \cdots & \cdots & \cdots \\
		 0 & 0 & 0 & \gamma & 1 \\
		 0 & 0 & 0 & 0 & \gamma
	  \end{array} \right) \]

\begin{lemma} \label{semisimple} An element $\alpha \in Mat(n,q)$ is semisimple if and only if
$\lambda_{\phi}(\alpha)_2'=0$ (i.e. all parts in all partitions in the
rational canonical form of $\alpha$ have size at most one).
\end{lemma}

\begin{proof}
	The explicit description of Jordan canonical form just given
implies that $\alpha$ is diagonalizable over $\bar{F_q}$ if and only if
there are no companion matrices $C(\phi_i^A)$ where $A \geq 2$.
\end{proof}

	Another tool we will use in applying the cycle index to semisimplicity
is Gordon's generalization of the Rogers-Ramanujan identities. It is worth
remarking that this gives what seems to be the first appearance of these
identities in finite group theory (see Fulman \cite{fulRogers} for more on
this). Lemma $\ref{Gordon}$ is a statement of Gordon's generalization of the
Rogers-Ramanujan identities. It is taken directly from page 111 of Andrews
$\cite{And}$.

\begin{lemma} \label{Gordon} For $1\leq i \leq k, k \geq 2$, and $|x|<1$

\[ \sum_{n_1,\cdots,n_{k-1} \geq 0} \frac{x^{N_1^2 + \cdots + N_{k-1}^2+
N_i + \cdots+ N_{k-1}}}{(x)_{n_1}\cdots(x)_{n_{k-1}}} = \prod_{r=1 \atop r
\neq 0, \pm i (mod \ 2k+1)}^{\infty} \frac{1}{1-x^r} \]

	where $N_j = n_j + \cdots n_{k-1}$.

\end{lemma}

	Lemma $\ref{Gordon}$ yields the following corollary.

\begin{cor} \label{part}

\[ \prod_{r=1}^{\infty} (1-\frac{1}{q^{rm_{\phi}}}) \sum_{\lambda: \lambda_1<k}
\frac{1}{c_{GL,\phi,q^{m_{\phi}}}(\lambda)} = \prod_{r=1 \atop r=0, \pm k (mod \
2k+1)}^{\infty} (1-\frac{1}{q^{m_{\phi}r}}) \]

\end{cor}

\begin{proof}
	Assume without loss of generality that $\phi=z-1$. By the second equality in
Theorem
$\ref{Rewrite}$, this probability is equal to:

\begin{eqnarray*}
\prod_{r=1}^{\infty} (1-\frac{1}{q^r}) \sum_{\lambda:
\lambda_1<k} \frac{1}{c_{GL,\phi,q^{m_{\phi}}}(\lambda)} & = & \prod_{r=1}^{\infty}
(1-\frac{1}{q^r}) \sum_{\lambda: \lambda_1<k} \frac{1}{q^{\sum_i
(\lambda_i')^2}
\prod_i (\frac{1}{q})_{m_i(\lambda)}}\\ & = & \prod_{r=1}^{\infty}
(1-\frac{1}{q^r}) \sum_{m_1(\lambda),\cdots,m_{k-1}(\lambda) \geq 0}
\frac{1}{q^{\sum_i (\lambda_i')^2} \prod_i (\frac{1}{q})_{m_i(\lambda)}}
\end{eqnarray*}

	The result now follows from Lemma $\ref{Gordon}$ by setting
$n_i=m_i(\lambda)$, $i=k$, and $x=\frac{1}{q}$.
\end{proof}

	With these tools in hand, the main results of this application can now be
obtained.

\begin{theorem} \label{Matss} The $n \rightarrow \infty$ limiting
probability that an element of $Mat(n,q)$ is semisimple is:

\[ \prod_{r=1 \atop r=0,\pm 2 (mod \ 5)}^{\infty} (1-\frac{1}{q^{r-1}}) \]

\end{theorem}

\begin{proof}
	By the cycle index for $Mat(n,q)$, Lemmas $\ref{bign}$, $\ref{allpoly}$,
\ref{product}, $\ref{semisimple}$ and Corollary $\ref{part}$, the
$n \rightarrow \infty$ limit of the chance that an element of $Mat(n,q)$ is
semisimple is:

\begin{eqnarray*}
& & lim_{n \rightarrow \infty} \frac{|GL(n,q)|}{q^{n^2}} [u^n]
\prod_{\phi} \sum_{\lambda:
\lambda_1<2} \frac{u^{m_{\phi}|\lambda|}}{c_{GL,\phi,q^{m_{\phi}}}(\lambda)}\\ 
& = & lim_{n \rightarrow \infty} \frac{|GL(n,q)|}{q^{n^2}} [u^n] \frac{1}{1-u}
\frac{1}{\prod_{r=1}^{\infty} (1-\frac{u}{q^r})}
\prod_{\phi} [\prod_{r=1}^{\infty} (1-\frac{u^{m_{\phi}}}{q^{rm_{\phi}}})
\sum_{\lambda:
\lambda_1<2} \frac{u^{m_{\phi}|\lambda|}}{c_{GL,\phi,q^{m_{\phi}}}(\lambda)}]\\ 
& = & \prod_{\phi} [\prod_{r=1}^{\infty} (1-\frac{1}{q^{rm_{\phi}}})
\sum_{\lambda:
\lambda_1<2} \frac{1}{c_{GL,\phi,q^{m_{\phi}}}(\lambda)}]\\ 
& = & \prod_{r=1 \atop r=0,\pm 2 (mod
\ 5)}^{\infty} (1-\frac{1}{q^{r-1}})
\end{eqnarray*}

\end{proof}

\begin{theorem} \label{ss} The $n \rightarrow \infty$ limiting probability that
an element of $GL(n,q)$ is semisimple is:

\[ \prod_{r=1 \atop r=0,\pm 2 (mod \ 5)}^{\infty}
\frac{(1-\frac{1}{q^{r-1}})}{(1-\frac{1}{q^r})} \]

\end{theorem}

\begin{proof}
	Arguing as in Theorem $\ref{Matss}$ the chance is:

\begin{eqnarray*}
& & lim_{n \rightarrow \infty} [u^n] \prod_{\phi \neq z}
\sum_{\lambda:\lambda_1<2}
\frac{u^{m_{\phi}|\lambda|}}{c_{GL,\phi,q^{m_{\phi}}}(\lambda)}\\
& = & lim_{n \rightarrow \infty} [u^n] \frac{1}{1-u} \prod_{\phi \neq z}
\prod_{r=1}^{\infty} (1-\frac{u^{m_{\phi}}}{q^{rm_{\phi}}})
\sum_{\lambda:
\lambda_1<2}\frac{u^{m_{\phi}|\lambda|}}{c_{GL,\phi,q^{m_{\phi}}}(\lambda)}\\
& = & \prod_{\phi \neq z} \prod_{r=1}^{\infty}
(1-\frac{1}{q^{rm_{\phi}}}) \sum_{\lambda: \lambda_2'=0}
\frac{1}{c_{GL,\phi,q^{m_{\phi}}}(\lambda)}]\\ 
& = & \prod_{\phi \neq z} \prod_{r=1 \atop r=0,\pm 2 (mod \ 5)}^{\infty}
(1-\frac{1}{q^{r}})\\
& = &  \prod_{r=1 \atop r=0,\pm 2 (mod \ 5)}^{\infty}
(\frac{1}{1-\frac{1}{q^{r}}}) \prod_{\phi} \prod_{r=1 \atop r=0,\pm 2 (mod \
5)}^{\infty} (1-\frac{1}{q^{r}})\\
& = & \prod_{r=1 \atop r=0,\pm 2 (mod \
5)}^{\infty} \frac{(1-\frac{1}{q^{r-1}})}{(1-\frac{1}{q^r})}
\end{eqnarray*}

\end{proof}

\begin{center}
{\bf Application 2: Regularity and Regular-Semisimplicity}
\end{center}

	An element of an algebraic group $G$ (usually taken to be connected
and reductive) is called regular if the dimension of its centralizer is as
small as possible (this minimal dimension turns out to be equal to the rank
of $G$, i.e. the dimension of a maximal torus of $G$). Regular elements are
very important in the representation theory of finite groups of Lie
type. This application will give formulas for the $n \rightarrow \infty$ limit of
the chance that an element of $GL(n,q)$ or $Mat(n,q)$ is regular or
regular-semisimple.

	The first step is to translate the conditions of regularity and regular
semisimplicity into conditions on the parititions $\lambda_{\phi}$. Let
$\bar{F_q}$ denote an algebraic closure of the field of $q$ elements.

\begin{lemma} \label{reg} An element $\alpha \in GL(n,q) \subset
GL(n,\bar{F_q})$ is regular if and only if all
$\lambda_{\phi}(\alpha)$ have at most one part.  \end{lemma}

\begin{proof}
	Let $\beta$ be an eigenvalue of $\alpha$ over $\bar{F_q}$ and
let $V_{\beta}$ be the eigenspace associated to $\beta$. The dimension
of $V_{\beta}$ is $|\lambda_{\phi}(\alpha)|$. Let $\alpha
|_{V_{\beta}}$ be the restriction of $\alpha$ to $V_{\beta}$. It is
not hard to see that:

\[ C_{GL(n,\bar{F_q})} (\alpha) = \prod_{\phi} \prod_{\beta \ root \ of \ 
\phi} C_{GL(|\lambda_{\phi}(\alpha)|,\bar{F_q})} (\alpha |_{V_{\beta}}) \]

	One can prove from Jordan canonical form, described before Lemma
\ref{semisimple} in the previous application, that
$C_{GL(|\lambda_{\phi}(\alpha)|,\bar{F_q})} (\alpha |_{V_{\beta}})$
has dimension $\sum_i (\lambda_{\phi,i}'(\alpha))^2$. Thus the
centralizer of $\alpha \in GL(n,\bar{F_q})$ has dimension $\sum_{\phi}
m_{\phi} \sum_i (\lambda_{\phi,i}'(\alpha))^2$. Given the value
$|\lambda_{\phi}(\alpha)|$, Lagrange multipliers show that $\sum_i
(\lambda_{\phi,i}'(\alpha))^2$ is minimized when
$\lambda_{\phi}(\alpha)$ has at most 1 part. The result follows since
$\sum_{\phi} m_{\phi} \sum_i \lambda_{\phi,i}'(\alpha)=n$.
\end{proof}

\begin{cor} \label{rss} An element of $GL(n,q)$ is regular-semisimple if and only
if
$|\lambda_{\phi}|=0,1$ for all $\phi$.
\end{cor}

\begin{proof}
	This is immediate from Lemmas \ref{semisimple} and \ref{reg}.
\end{proof}

	It is worth remarking that the condition of Lemma $\ref{reg}$,
and hence the condition of regularity in $GL(n,q)$, is equivalent to the condition
that the minimum polynomial of $\alpha$ is equal to the characteristic
polynomial of $\alpha$.

	Lemma $\ref{reg}$ leads us to call an element $\alpha \in
Mat(d,q)$ regular if all $\lambda_{\phi}(\alpha)$ have at most one
part. Neumann and Praeger $\cite{Nu1},\cite{Nu2}$ studied the chance
that a matrix is regular or regular-semisimple (they called these
conditions "cyclic" and "separable"). They were interested in these
probabilities because they give a way to test random number generators
and computer algorithms for generating random elements from a finite
group (see Chapter 1 of Fulman \cite{fulthesis} for further motivation). Volkmann
and Fleischmann
$\cite{Fle}$ and Lehrer
$\cite{Leh}$ studied this problem as well.

	Some theorems of Neumann and Praeger are:

\begin{enumerate}

\item For $n \geq 2$, the chance that an $n*n$ matrix is not regular
semi-simple is at least $q^{-1}-q^{-2}-q^{-3}$ and at most $q^{-1}+O(q^{-2})$.

\item For $n \geq 2$, the chance that an $n*n$ matrix is not regular
is at least $\frac{1}{q^2(q+1)}$ and at most $\frac{1}{(q^2-1)(q-1)}$.

\end{enumerate}

	In the next four theorems, the cycle index machinery is used
to find $n \rightarrow \infty$ formulas for the chance that a $n*n$
matrix or an element of $GL(n,q)$ is regular or
regular-semisimple. These are good examples of results which seem hard
to prove by other methods. Throughout, Lemma $\ref{bign}$ is used
freely.

\begin{theorem} \label{regssMat} The $n \rightarrow \infty$ chance that an
$n*n$ matrix in $Mat(n,q)$ is regular-semisimple is equal to:

\[ \prod_{r=1}^{\infty} (1-\frac{1}{q^r}) \]

\end{theorem}

\begin{proof}
	Corollary \ref{rss} says that an element of $Mat(n,q)$ is regular semisimple iff
all $\lambda_{\phi}$ have size at most 1. So the cycle index for $Mat(n,q)$ and
Lemma
$\ref{allpoly}$ imply that the probability of regular-semisimplicity
is:

\begin{eqnarray*}
&  &  \lim_{n \rightarrow \infty}
\frac{|GL(n,q)|}{q^{n^2}} [u^n] \prod_{\phi} (1+\frac{u^{m_{\phi}}}{q^{m_{\phi}}-1})\\
& = &  \lim_{n \rightarrow \infty} \frac{|GL(n,q)|}{q^{n^2}} [u^n]
\frac{\prod_{\phi} (1+\frac{u^{m_{\phi}}}{q^{m_{\phi}}-1}) (1-\frac{u^{m_{\phi}}}{q^{m_{\phi}}})}{1-u}\\
& = & \prod_{r=1}^{\infty} (1-\frac{1}{q^r}) \prod_{\phi}  (1+\frac{1}{q^{m_{\phi}}-1}) (1-\frac{1}{q^{m_{\phi}}})\\
& = & \prod_{r=1}^{\infty} (1-\frac{1}{q^r})
\end{eqnarray*}

\end{proof}

\begin{theorem} \label{regMat} The $n \rightarrow \infty$ chance that an
$n*n$ matrix in $Mat(n,q)$ is regular is equal to:

\[ (1-\frac{1}{q^5}) \prod_{r=3}^{\infty} (1-\frac{1}{q^r}) \]

\end{theorem}

\begin{proof}
	Lemma $\ref{reg}$ shows that regularity is equivalent to all
$\lambda_{\phi}$ having at most 1 part. The cycle index for $Mat(n,q)$
and Lemma $\ref{allpoly}$ give that the probability is:

\begin{eqnarray*}
&  & \lim_{n \rightarrow \infty}
\frac{|GL(n,q)|}{q^{n^2}} [u^n] \prod_{\phi} (1+ \sum_{j=1}^{\infty} \frac{u^{jm_{\phi}}} {q^{jm_{\phi}-m_{\phi}}(q^{m_{\phi}}-1)})\\
& = & \prod_{r=1}^{\infty} (1-\frac{1}{q^r}) \lim_{n \rightarrow \infty} [u^n]
\frac{\prod_{\phi} (1-\frac{u^{m_{\phi}}}{q^{m_{\phi}}}) (1+ \sum_{j=1}^{\infty} \frac{u^{jm_{\phi}}} {q^{jm_{\phi}-m_{\phi}}(q^{m_{\phi}}-1)})} {1-u}  \\
& = & \prod_{r=1}^{\infty} (1-\frac{1}{q^r}) \lim_{n \rightarrow \infty} [u^n]
\frac {\prod_{\phi} (1+ \frac{u^{m_{\phi}}}{q^{m_{\phi}}(q^{m_{\phi}}-1)})}{1-u}\\
& = & \prod_{r=1}^{\infty} (1-\frac{1}{q^r}) \prod_{\phi} (1+ \frac{1}{q^{m_{\phi}}(q^{m_{\phi}}-1)})\\
& = & \prod_{r=3}^{\infty} (1-\frac{1}{q^r}) \prod_{\phi} (1+
\frac{1}{q^{m_{\phi}}(q^{m_{\phi}}-1)})(1-\frac{1}{q^{2m_{\phi}}})(1-\frac{1}{q^{3m_{\phi}}})\\
& = &  \prod_{r=3}^{\infty} (1-\frac{1}{q^r}) \prod_{\phi}
(1-\frac{1}{q^{6m_{\phi}}})\\
& = &  (1-\frac{1}{q^5}) \prod_{r=3}^{\infty} (1-\frac{1}{q^r})
\end{eqnarray*}

\end{proof}

\begin{theorem} \label{regssGL} The $n \rightarrow \infty$ chance that an
element of $GL(n,q)$ is regular-semisimple is equal to:

\[ 1-\frac{1}{q} \]

\end{theorem}

\begin{proof}
	By the cycle index for $GL$, the probability is:

\begin{eqnarray*}
&  &  \lim_{n \rightarrow \infty} [u^n] \prod_{\phi \neq z} (1+\frac{u^{m_{\phi}}}{q^{m_{\phi}}-1})\\
& = &  \lim_{n \rightarrow \infty} [u^n] \frac{(1-\frac{u}{q})
\prod_{\phi \neq z}
[(1+\frac{u^{m_{\phi}}}{q^{m_{\phi}}-1})(1-\frac{u^{m_{\phi}}}{q^{m_{\phi}}})]}{1-u}\\
& = & 1-\frac{1}{q}
\end{eqnarray*}

\end{proof}

\begin{theorem} \label{regGL}
	The $n \rightarrow \infty$ chance that an element of $GL(n,q)$ is regular is equal to:

\[ \frac{1-\frac{1}{q^5}}{1+\frac{1}{q^3}} \]

\end{theorem}

\begin{proof}
	By the cycle index for $GL$, the probability is:

\begin{eqnarray*}	
&  & \lim_{n \rightarrow \infty} [u^n] \prod_{\phi \neq z} (1+\sum_{j=1}^{\infty} \frac{u^{jm_{\phi}}} {q^{jm_{\phi}-m_{\phi}}(q^{m_{\phi}}-1)})\\
& = & \lim_{n \rightarrow \infty} [u^n] \frac{(1-\frac{u}{q})  \prod_{\phi \neq z}
(1+\frac{u^{m_{\phi}}}{(q^{m_{\phi}}-1)(1-\frac{u^{m_{\phi}}}{q^{m_{\phi}}})})(1-\frac{u^{m_{\phi}}}{q^{m_{\phi}}})}{1-u}\\
& = & (1-\frac{1}{q}) \prod_{\phi \neq z}^{\infty}
(1+\frac{1}{q^{m_{\phi}}(q^{m_{\phi}}-1)})\\
& = &  \frac{1-\frac{1}{q}}{1+\frac{1}{q(q-1)}} \prod_{\phi}
(\frac{1-\frac{1}{q^{6m_{\phi}}}}{(1-\frac{1}{q^{2m_{\phi}}})(1-\frac{1}{q^{3m_{\phi}}})})\\
& = &  \frac{1-\frac{1}{q}}{1+\frac{1}{q(q-1)}}
\frac{1-\frac{1}{q^5}}{(1-\frac{1}{q})(1-\frac{1}{q^2})}\\
& = & \frac{1-\frac{1}{q^5}}{1+\frac{1}{q^3}}
\end{eqnarray*}

\end{proof}

	It is worth remarking, along the lines of Neumann and Praeger
$\cite{Nu1}$, that these results are intuitively reasonable. Namely,
Steinberg $\cite{Ste}$ proved that the set of non-regular elements in
an algebraic group has codimension $3$. Suppose this set to be a
non-singular (which it is not) high dimensional (which it is for large
$n$) variety. Then the chance of non-regularity would be about
$\frac{1}{q^3}$, so the chance of regularity would be about
$1-\frac{1}{q^3}$ which is consistent with Theorems $\ref{regMat}$ and
$\ref{regGL}$.

	Similarly, Neumann and Praeger $\cite{Nu1}$ noted that a
matrix is regular semisimple if and only if the discriminant of its
characteristic polynomial is non-zero. As this restriction is defined
by one equation, the chance of being regular semisimple should be
about $1-\frac{1}{q}$, which is consistent with Theorems
$\ref{regssMat}$ and $\ref{regssGL}$.

	Finally, the chance of being semi-simple is greater than the
chance of being regular semi-simple, which is consistent with Theorem
$\ref{ss}$.

\section{Cycle Indices of the Finite Unitary, Symplectic, and Orthogonal
Groups} \label{OTHER}

	This section develops cycle indices for the finite unitary, symplectic, and
orthogonal groups, and in that order. For the symplectic and orthogonal groups,
it will be assumed that the characteristic is not equal to 2. The main ingredient
is work of Wall \cite{Wal} on the conjugacy classes of these groups. The
rewritings of
$c_{GL,\phi,q^{m_{\phi}}}$ from Section \ref{GLCYC} will be used to relate these cycle
indices to the cycle index of the general linear groups.

\begin{center}
{\bf The Unitary Groups}
\end{center}

	The unitary group $U(n,q)$ (we allow characteristic 2) can be defined
as the subgroup of $GL(n,q^2)$ preserving a non-degenerate skew-linear form.
Recall that a skew-linear form on an
$n$ dimensional vector space $V$ over $F_{q^2}$ is a bilinear map
$<,>:V \times V \rightarrow F_{q^2}$ such that $<\vec{x},\vec{y}> =
<\vec{y},\vec{x}>^q$ (raising to the $q$th power is an involution in a
field of order $q^2$). One such form is given by $<\vec{x},\vec{y}> =
\sum_{i=1}^n x_i y_i^q$. It is known (page 7 of Carter $\cite{Ca1}$)
that any two non-degenerate skew-linear forms are equivalent, so that
$U(n,q)$ is unique up to isomorphism. The
order of $U(n,q)$ is $q^{{n \choose 2}} \prod_{i=1}^n (q^i - (-1)^i)$.

	To define a cycle index for the finite unitary groups, it is first necessary to
understand their conjugacy classes. Given a polynomial $\phi$ with coefficients in
$F_{q^2}$ and non vanishing constant term, define a polynomial $\tilde{\phi}$ by:

\[ \tilde{\phi} = \frac{z^{m_{\phi}} \phi^q(\frac{1}{z})}{[\phi(0)]^q} \]

	where $\phi^q$ raises each coefficient of $\phi$ to the $q$th
power. Writing this out, a polynomial $\phi(z)=z^{m_{\phi}} +
\alpha_{m_{\phi}-1} z^{m_{\phi}-1} + \cdots + \alpha_1 z + \alpha_0$
with $\alpha_0 \neq 0$ is sent to $\tilde{\phi}(z)= z^{m_{\phi}} +
(\frac{\alpha_1}{\alpha_0})^q z^{m_{\phi}-1}+ \cdots +
(\frac{\alpha_{m_{\phi}-1}} {\alpha_0})^qz + (\frac{1}{\alpha_0})^q$.

	Wall $\cite{Wal}$ proves that the conjugacy classes of the
unitary group correspond to the following combinatorial data. As was
the case with $GL(n,q^2)$, an element $\alpha \in U(n,q)$ associates
to each monic, non-constant, irreducible polynomial $\phi$ over
$F_{q^2}$ a partition $\lambda_{\phi}$ of some non-negative integer
$|\lambda_{\phi}|$ by means of rational canonical form. The
restrictions necessary for the data $\lambda_{\phi}$ to represent a
conjugacy class are:

\begin{enumerate}
\item $|\lambda_z|=0$
\item $\lambda_{\phi}=\lambda_{\tilde{\phi}}$
\item $\sum_{\phi} |\lambda_{\phi}|m_{\phi}=n$
\end{enumerate}

	Recall that $m_i(\lambda)$ denotes the number of parts of
$\lambda$ of size $i$. Wall computed the size of a conjugacy class
corresponding to the data $\lambda_{\phi}$ as:

\[ \frac{|U(n,q)|}{ \prod_{\phi} B(\phi)} \]

	where

\begin{eqnarray*}
A(\phi^i) & = & |U(m_i(\lambda_{\phi}),q^{m_{\lambda}})| \ if \
\phi=\tilde{\phi}\\
& = & |GL(m_i(\lambda_{\phi}),q^{2m_{\lambda}})|^{\frac{1}{2}} \ if \ \phi \neq
\tilde{\phi}\\
\end{eqnarray*} 

	and

\[ B(\phi) = q^{2m_{\phi}[\sum_{h<i}
hm_h(\lambda_{\phi})m_i(\lambda_{\phi}) + \frac{1}{2} \sum_i
(i-1)m_i(\lambda_{\phi})^2]} \prod_i A(\phi^i). \]

	As an example of this formula, consider the set of unitary
transvections, namely determinant 1 elements of $U(n,q)$ whose
pointwise fixed space is $n-1$ dimensional. Fulman \cite{fulthesis}
counts unitary transvections directly, showing that the number of them is:

\[ \frac{(q^n-(-1)^n)(q^{n-1}-(-1)^{n-1})}{q+1} \]

	This checks with the class size formula. The fact that the dimension of the
fixed space of $\alpha \in U(n,q)$ is the number of parts of
$\lambda_{z-1}(\alpha)$ implies that $\alpha \in U(n,q)$ is a unitary
transvection exactly when $\lambda_{z-1}(\alpha)=(2,1^{n-2})$ and
$|\lambda_{\phi}(\alpha)|=0$ for $\phi \neq z-1$. Thus Wall's work implies that
the unitary transvections form a single conjugacy class of size:

\[ \frac{|U(n,q)|}{q^{2n-3} |U(n-2,q)| |U(1,q)|} =
\frac{(q^n-(-1)^n)(q^{n-1}-(-1)^{n-1})}{q+1}. \]

	In analogy with the general linear groups, define a cycle index for the unitary
groups by:

	\[ 1 + \sum_{n=1}^{\infty} \frac{u^n}{|U(n,q)|} \sum_{\alpha \in U(n,q)}
\prod_{\phi \neq z} x_{\phi,\lambda_{\phi}(\alpha)} \]

	Theorem $\ref{cycleU}$ will give a factorization theorem for this
cycle index in terms of quantities for $GL$ (so that some of the
results obtained for the general linear groups can be carried
over for free). For this and future use in the applications of Section
\ref{APPLICATIONS}, it is desirable to count the number of $\phi$ of a given
degree invariant under $\tilde{}$.

\begin{lemma} \label{pro} $\widetilde{\phi_1 \phi_2} = \tilde{\phi_1}
\tilde{\phi_2}$.
\end{lemma}

\begin{proof}
	From the definition of the involution $\tilde{}$, the lemma
reduces to the observation that $(\phi_1 \phi_2)^q =
(\phi_1^q)(\phi_2^q)$, where $^q$ is the map which raises each
coefficient of a polynomial to the $q$th power.  \end{proof}

	Let $\tilde{I}_{m,q^2}$ be the number of monic, irreducible
polynomials $\phi$ of degree $m$ over $F_{q^2}$ such that
$\phi=\tilde{\phi}$.

\begin{theorem} \label{Count1} $\tilde{I}_{m,q^2}=0$ if $m$ is even and
$\tilde{I}_{m,q^2}=\frac{1}{m} \sum_{d|m} \mu(d) (q^{\frac{m}{d}}+1)$
if $m$ is odd.
\end{theorem}

\begin{proof}
	Let $M_m$ be the number of monic degree $m$ polynomials (not
necessarily irreducible) and let $\tilde{M}_m$ be the number of monic
degree $m$ polynomials $\phi$ (not necessarily irreducible) such that
$\phi(0) \neq 0$ and $\phi=\tilde{\phi}$. Define $A(t)= 1 +
\sum_{m=1}^{\infty} M_m t^m$ and $B(t)= 1 + \sum_{m=1}^{\infty}
\tilde{M}_m t^m$. Note that $A(t)=\frac{1}{1-q^2t}$ because
$M_m=q^{2m}$. Wall $\cite{Wal}$ observes that $B(t)=\frac{1+t}{1-qt}$
(this follows from the fact that $\tilde{M}_m=q^m+q^{m-1}$, which is
clear from the explicit description of the definition of
$\tilde{\phi}$ given above).

	The fact that the involution $\tilde{}$ preserves degree gives
the following equation (where as usual all polynomials in the products
are irreducible):

	\[ A(t) = \frac{1}{1-t} \prod_{\phi \neq z,\phi=\tilde{\phi}}
(1+\sum_{n=1}^{\infty} t^{n m_{\phi}}) \prod_{\{\phi,\tilde{\phi}\},
\phi \neq \tilde{\phi}} (1+\sum_{n=1}^{\infty} t^{n m_{\phi}})^2 \]

	Lemma $\ref{pro}$ implies that a polynomial invariant under
$\tilde{}$ is a product of terms $\phi$ where $\phi=\tilde{\phi}$ and $\phi
\tilde{\phi}$ where $\phi \neq \tilde{\phi}$. This gives the equation:

	\[ B(t) = \prod_{\phi \neq z,\phi=\tilde{\phi}}
(1+\sum_{n=1}^{\infty} t^{n m_{\phi}}) \prod_{\{\phi,\tilde{\phi}\},
\phi \neq \tilde{\phi}} (1+\sum_{n=1}^{\infty} t^{2nm_{\phi}}) \]

	These equations give:

\begin{eqnarray*}
\frac{B(t)^2}{A(t^2)} & = & (1-t^2) \frac{\prod_{m=0}^{\infty}
(1+\sum_{n=1}^{\infty} t^{mn})^{2\tilde{I}_{m,q^2}}}{\prod_{m=0}^{\infty}
(1+\sum_{n=1}^{\infty} t^{2mn})^{\tilde{I}_{m,q^2}}}\\
& = & (1-t^2) \prod_{m=0}^{\infty} \frac{(1-t^{2m})^{\tilde{I}_{m,q^2}}}{(1-t^m)^{2\tilde{I}_{m,q^2}}}\\
& = & (1-t^2) \prod_{m=0}^{\infty} (\frac{1+t^m}{1-t^m})^{\tilde{I}_{m,q^2}}
\end{eqnarray*}

	Combining this with the explicit expressions for $A(t)$ and $B(t)$
given above shows that:

	\[ \prod_{m=0}^{\infty} (\frac{1+t^m}{1-t^m})^{\tilde{I}_{m,q^2}} =
(\frac{1+t}{1-t})(\frac{1+qt}{1-qt}) \]

	Take logarithms of both sides of this equation, using the
expansions $log(1+x)=x-\frac{x^2}{2}+\frac{x^3}{3}+\cdots$ and
$log(1-x)=-x-\frac{x^2}{2}-\frac{x^3}{3}+\cdots$.

	The left-hand side becomes:

	\[ \sum_{m=0}^{\infty}
2 \tilde{I}_{m,q^2}(t^m+\frac{t^{3m}}{3}+\frac{t^{5m}}{5}+\cdots) \]

	The right-hand becomes:

	\[ \sum_{m \ odd} 2(\frac{1+q^m}{m}) t^m \]

	Comparing coefficients of $t^m$ shows that
$\tilde{I}_{m,q^2}=0$ for $m$ even and that $\sum_{d|m} 2
\tilde{I}_{d,q^2} \frac{d}{m} = 2 (\frac{1+q^m}{m})$ for $m$
odd. Moebius inversion proves that $\tilde{I}_{m,q^2}=\frac{1}{m}
\sum_{d|m} \mu(d) (1+q^{\frac{m}{d}})$ if $m$ is odd.
\end{proof}

	Next, it will be proved that the cycle index of the unitary
groups factors. (One can prove a factorization theorem without Theorem
$\ref{Count1}$, but Theorem $\ref{Count1}$ is necessary to get an
expression in terms of quantities related to $GL$). The quantities
$c_{GL,\phi,q^{m_{\phi}}}$ and their various rewritings were considered in
Section $\ref{GLCYC}$.

\begin{theorem} \label{cycleU}

\begin{eqnarray*}
1+\sum_{n=1}^{\infty} \frac{u^n}{|U(n,q)|} \sum_{\alpha \in U(n,q)}
\prod_{\phi \neq z} x_{\phi,\lambda_{\phi}(\alpha)} & = & \prod_{\phi \neq z,
\phi=\tilde{\phi}} [\sum_{\lambda} x_{\phi,\lambda}
\frac{(-u)^{|\lambda|m_{\phi}}} {c_{GL,z-1,-(q^{m_{\phi}})}(\lambda)}]\\
& & \prod_{\{\phi,\tilde{\phi}\}, \phi \neq \tilde{\phi}} [\sum_{\lambda}
x_{\phi,\lambda} x_{\tilde{\phi},\lambda}
\frac{u^{2|\lambda|m_{\phi}}}{c_{GL,z-1,q^{2m_{\phi}}}(\lambda)}]
\end{eqnarray*}

\end{theorem}

\begin{proof}
	The theorem follows from Wall's description and formula for
conjugacy class sizes in the unitary group, provided that one can prove that for
all
$\phi=\tilde{\phi}$,

\begin{eqnarray*}
& & \frac{u^{|\lambda_{\phi}|}}{q^{2m_{\phi}[\sum_{h<i}
hm_h(\lambda_{\phi})m_i(\lambda_{\phi}) + \frac{1}{2} \sum_i
(i-1)m_i(\lambda_{\phi})^2]} \prod_i |U(m_i(\lambda_{\phi}),q^{m_{\phi}})|}\\
& = & \frac{(-u)^{|\lambda_{\phi}|}} { (-q)^{2m_{\phi}[\sum_{h<i}
hm_h(\lambda_{\phi})m_i(\lambda_{\phi}) + \frac{1}{2} \sum_i
(i-1)m_i(\lambda_{\phi})^2]} \prod_i
|GL(m_i(\lambda_{\phi}),(-q)^{m_{\phi}})|}
\end{eqnarray*}

	The formulas for $|GL(n,q)|$ and $|U(n,q)|$ show that
$|GL(n,-q)| = (-1)^n |U(n,q)|$.

	The proof of the desired equation boils down to keeping track
of powers of $-1$ and using the fact from Theorem $\ref{Count1}$ that
if $\phi= \tilde{\phi}$, then $\phi$ has odd degree. With a little
more detail,

\begin{eqnarray*}
|\lambda_{\phi}|+ m_{\phi} [\sum_i
(i-1)m_i(\lambda_{\phi})^2] + m_{\phi}[\sum_i m_i(\lambda_{\phi})] & = & |\lambda_{\phi}|+ \sum_i
(i-1)m_i(\lambda_{\phi})^2 + m_i(\lambda_{\phi}) \ (mod \ 2)\\
& = &  |\lambda_{\phi}| + \sum_i im_i(\lambda_{\phi}) \ (mod \ 2)\\
& = & 2|\lambda_{\phi}| \ (mod \ 2)\\
& = & 0 \ (mod \ 2)
\end{eqnarray*}

\end{proof}

\begin{center}
{\bf The Symplectic Groups}
\end{center}

	This paper assumes for simplicity that the characteristic of $F_q$
is not equal to 2. The symplectic group $Sp(2n,q)$ can be defined as
the subgroup of $GL(2n,q)$ preserving a non-degenerate alternating
form on $F_q$. Recall that an alternating form on a $2n$ dimensional
vector space $V$ over $F_q$ is a bilinear map $<,>:V \times V
\rightarrow F_q$ such that $<\vec{x},\vec{y}>=-<\vec{y},\vec{x}>$
(alternating forms do not exist in odd dimension). One such form is
given by $<\vec{x},\vec{y}> = \sum_{i=1}^n (x_{2i-1}y_{2i} -
x_{2i}y_{2i-1})$. As is explained in Chapter 1 of Carter $\cite{Ca1}$,
there is only one such form up to equivalence, so $Sp(2n,q)$ is unique
up to isomorphism. The order of $Sp(2n,q)$ is $q^{n^2} \prod_{i=1}^n (q^{2i}-1)$.
	
	With the aim of finding a cycle index for the symplectic groups, let us
first understand their conjugacy classes. Given a a polynomial
$\phi$ with coefficients in
$F_q$ and non vanishing constant term, define a polynomial $\bar{\phi}$ by:

\[ \bar{\phi} = \frac{z^{m_{\phi}} \phi^q(\frac{1}{z})}{[\phi(0)]^q} \]

	where $\phi^q$ raises each coefficient of $\phi$ to the $q$th
power. Explicitly, a polynomial $\phi(z)=z^{m_{\phi}} +
\alpha_{m_{\phi}-1} z^{m_{\phi}-1} + \cdots + \alpha_1 z + \alpha_0$
with $\alpha_0 \neq 0$ is sent to $\bar{\phi}(z)= z^{m_{\phi}} +
(\frac{\alpha_1}{\alpha_0})^q z^{m_{\phi}-1}+ \cdots +
(\frac{\alpha_{m_{\phi}-1}} {\alpha_0})^qz +
(\frac{1}{\alpha_0})^q$. The notation $\bar{\phi}$ breaks from Wall
$\cite{Wal}$, in which $\tilde{\phi}$ was used, but these maps are
different. Namely $\tilde{}$ is defined on polynomials with
coefficients in $F_q$, but $\bar{}$ is defined on polynomials with
coefficients in $F_{q^2}$. The distinction between the maps $\tilde{}$
and $\bar{}$ will be evident in the different statements of Theorems
$\ref{Count1}$ and $\ref{Count2}$.

	Wall $\cite{Wal}$ showed that a conjugacy class of $Sp(2n,q)$
corresponds to the following data. To each monic, non-constant,
irreducible polynomial $\phi \neq z \pm 1$ associate a partition
$\lambda_{\phi}$ of some non-negative integer $|\lambda_{\phi}|$. To
$\phi$ equal to $z-1$ or $z+1$ associate a symplectic signed partition
$\lambda_{\phi}^{\pm}$, by which is meant a partition of some natural
number $|\lambda_{\phi}^{\pm}|$ such that the odd parts have even
multiplicity, together with a choice of sign for the set of parts of
size $i$ for each even $i>0$.

\begin{center}
Example of a Symplectic Signed Partition
\end{center}
	
\[ \begin{array}{c c c c c c c c}
	&  . & . & . & . & .  \\
	&  . & . & . & . & .  \\
	+  & . & . & . & . &  \\
	&  . & . & . &  &  \\
	&  . & . & . &  &   \\
	-  & . & . & & &   \\
	 & . & . & & &  
	  \end{array}  \]

	Here the $+$ corresponds to the parts of size 4 and the $-$
corresponds to the parts of size 2. This data represents a conjugacy
class of $Sp(2n,q)$ if and only if:

\begin{enumerate}

\item $|\lambda_{z}|=0$
\item $\lambda_{\phi}=\lambda_{\bar{\phi}}$
\item $\sum_{\phi=z \pm 1} |\lambda_{\phi}^{\pm}| + \sum_{\phi \neq z \pm 1} |\lambda_{\phi}| m_{\phi}=2n$

\end{enumerate}

	Wall computed the size of a conjugacy class corresponding to this
data as:

\[ \frac{|Sp(2n,q)|}{\prod_{\phi} B(\phi)} \]

	where

\begin{eqnarray*}
B(\phi) = q^{[\sum_{h<i}
hm_h(\lambda_{\phi}^{\pm})m_i(\lambda_{\phi}^{\pm}) + \frac{1}{2} \sum_i
(i-1)m_i(\lambda_{\phi}^{\pm})^2]} \prod_i A(\phi^{\pm,i}) \ if \ \phi=z
\pm 1\\
B(\phi) = q^{m_{\phi}[\sum_{h<i}
hm_h(\lambda_{\phi})m_i(\lambda_{\phi}) + \frac{1}{2} \sum_i
(i-1)m_i(\lambda_{\phi})^2]} \prod_i A(\phi^i) \ if \ \phi \neq z \pm 1
\end{eqnarray*}
  
	and

\begin{eqnarray*}
A(\phi^{\pm,i}) & = & |Sp(m_i(\lambda_{\phi}^{\pm}),q)| \ if
\ i=1 \ (mod \ 2)\\
& = & q^{\frac{m_i(\lambda_{\phi}^{\pm})}{2}}
|O(m_i(\lambda_{\phi}^{\pm}),q)| \ if \  i=0 \ (mod \ 2)\\
A(\phi^i) & = & |U(m_i(\lambda_{\phi}),q^{\frac{m_{\lambda}}{2}})| \ if \
\phi=\bar{\phi}\\
& = & |GL(m_i(\lambda_{\phi}),q^{m_{\lambda}})|^{\frac{1}{2}} \ if \ \phi \neq
\bar{\phi}.
\end{eqnarray*}

	Here $O(m_i(\lambda_{\phi}),q)$ is the orthogonal group with
the same sign as the sign associated to the parts of size $i$ (see the upcoming
treatment of the orthogonal groups in this section for more background
on them). The quantity
$B(\phi)$ will also be denoted by $c_{Sp,\phi,q^{m_{\phi}}}(\lambda^{\pm})$.

	As an example, consider the set of symplectic transvections,
i.e. determinant 1 elements of $Sp(2n,q)$ whose pointwise fixed space
is $n-1$ dimensional. Fulman \cite{fulthesis} uses direct counting arguments to
show that there are a total of $q^{2n}-1$ symplectic transvections.

	These symplectic transvections split into two conjugacy
classes of size $\frac{q^{2n}-1}{2}$. Using the fact that the dimension of the
fixed space of
$\alpha$ is the number of parts of $\lambda^{\pm}_{z-1}(\alpha)$, one concludes
that an element of $Sp(2n,q)$ is a symplectic transvection if
and only if $\lambda_{z-1}(\alpha)$ is $(+2,1^{n-2})$ or
$(-2,1^{n-2})$ and $|\lambda_{\phi}|=0$ for $\phi \neq z-1$. By Wall's
class size formula, the sizes of these conjugacy classes are:

\[ \frac{|Sp(2n,q)|}{q^{2n-\frac{3}{2}} |Sp(2n-2,q)| q^{\frac{1}{2}}
|O^{\pm}(1,q)|} = \frac{q^{2n}-1}{2} \]

	for both conjugacy classes. This confirms that there are
$q^{2n}-1$ symplectic transvections.

	As with the general linear and unitary groups define a cycle index for the
symplectic groups by

\[ 1 + \sum_{n=1}^{\infty} \frac{u^{2n}}{|Sp(2n,q)|} \sum_{\alpha \in
Sp(2n,q)}
\prod_{\phi = z \pm 1} x_{\phi,\lambda^{\pm}_{\phi}(\alpha)} \prod_{\phi
\neq z, z \pm 1} x_{\phi,\lambda_{\phi}(\alpha)} \]

	Theorem $\ref{CycleSym}$ will prove that this cycle index
factors. For this Theorem $\ref{Count2}$, which counts polynomials
invariant under the involution $\bar{}$, is essential. Let
$\bar{I}_{m,q}$ be the number of monic irreducible polynomials $\phi$
of degree $m$ with coefficients in $F_q$ such that $\phi=\bar{\phi}$.

\begin{lemma} \label{product2} $\overline{\phi_1 \phi_2} =
\bar{\phi_1} \bar{\phi_2}$.
\end{lemma}

\begin{proof}
	From the definition of the involution $\bar{}$, the lemma
reduces to the observation that $(\phi_1 \phi_2)^q =
(\phi_1^q)(\phi_2^q)$, where $^q$ is the map which raises each
coefficient of a polynomial to the $q$th power.
\end{proof}
 
\begin{theorem} \label{Count2}

\begin{enumerate}

\item $\bar{I}_{1,q}=2$ and the two degree 1 polynomials such
that $\phi=\bar{\phi}$ are $z \pm 1$.

\item If $m \neq 1$ is odd, then $\bar{I}_{m,q}=0$.

\item If $m=2^r m_0$ is even, with $m_0$ odd, then $\bar{I}_{m,q}=\frac{1}{m}
\sum_{d|m_0} \mu(d) (q^{\frac{m}{2d}}-1)$.

\end{enumerate}
\end{theorem}

\begin{proof}
	The method of proof is essentially the same as that used for
the unitary groups in Theorem $\ref{Count1}$. Let $M_m$ be the number
of monic degree $m$ polynomials (not necessarily irreducible) over
$F_q$ and let $\bar{M}_m$ be the number of monic degree $m$
polynomials $\phi$ (not necessarily irreducible) over $F_q$ such that
$\phi(0) \neq 0$ and $\phi=\bar{\phi}$. Define
$A(t)= 1 + \sum_{m=1}^{\infty} M_m t^m$ and $B(t)= 1 + \sum_{m=1}^{\infty}
\bar{M}_m t^m$. Note that $A(t)=\frac{1}{1-qt}$ because $M_m=q^m$. On
page 37 of Wall $\cite{Wal}$ it is noted that
$B(t)=\frac{(1+t)^2}{1-qt^2}$ (this follows from the explicit
description of the definition of $\bar{\phi}$.

	Arguing exactly as for the unitary group in Theorem
$\ref{Count1}$ gives:

\[ \frac{B(t)^2}{A(t^2)}  =  (1-t^2) \prod_{m=0}^{\infty}
(\frac{1+t^m}{1-t^m})^{\bar{I}_{m,q}} \]

	Combining this with the explicit expressions for $A(t)$ and $B(t)$
given above yields:

	\[ \prod_{m=0}^{\infty} (\frac{1+t^m}{1-t^m})^{\bar{I}_{m,q}} =
\frac{(1+t)^3}{(1-t)(1-qt^2)} \]

	Take logarithms of both sides of this equation, using the
expansions $log(1+x)=x-\frac{x^2}{2}+\frac{x^3}{3}+\cdots$ and
$log(1-x)=-x-\frac{x^2}{2}-\frac{x^3}{3}+\cdots$.

	The left-hand side becomes:

	\[ \sum_{m=0}^{\infty}
2\bar{I}_{m,q}(t^m+\frac{t^{3m}}{3}+\frac{t^{5m}}{5}+\cdots) \]

	The right-hand side becomes:

	\[ 4 \sum_{m \ odd} \frac{t^m}{m} + 2 \sum_{m \ even}
\frac{t^m}{m} (q^{\frac{m}{2}}-1) \]

	Comparing coefficients of $t$ shows that
$\bar{I}_{1,q}=2$. Since $z-1$ and $z+1$ satisfy $\phi=\bar{\phi}$,
these are the two degree 1 polynomials satisfying
$\phi=\bar{\phi}$. The $\bar{I}_{m,q}$ are all non-negative and the
odd degree terms on the right-hand side have been accounted for. Thus
$\bar{I}_{m,q}=0$ if $m \neq 1$ is odd.

	Now suppose that $m$ is even and write $m=2^rm_0$ where $m_0$
is odd. The coefficient of $t^m$ on the right-hand side is
$\frac{2}{m} (q^{\frac{m}{2}}-1)$. The coefficient of $t^m$ in the
left-hand side is $\sum_{k|m_0} 2 \frac{k}{m_0}
\bar{I}_{2^rk,q}$. This gives the relation:

\[ \sum_{k|m_0} k\bar{I}_{2^rk,q} = \frac{q^{2^{r-1}m_0}-1}{2^r} \]

	 It is straightforward to check that on the lattice of odd
integers with divisibility as the inclusion relation, Moebius
inversion holds in the sense that if $F(n)=\sum_{d|n} f(d)$ for all
odd $n$, then $f(n) = \sum_{d|n} \mu(d) F(\frac{n}{d})$ for all odd
$n$. Fix $r$ and define functions on the lattice of odd integers
by $F_r(n)= \frac{q^{2^{r-1} n}-1} {2^r}$ and $f_r(n)=n
\bar{I}_{2^rn,q}$. The theorem follows by Moebius inversion.
\end{proof}

	It can now be seen that the cycle index for the symplectic
groups factors.

\begin{theorem} \label{CycleSym} 

\begin{eqnarray*}
& & 1 + \sum_{n=1}^{\infty} \frac{u^{2n}}{|Sp(2n,q)|} \sum_{\alpha
\in Sp(2n,q)}
\prod_{\phi = z \pm 1} x_{\phi,\lambda^{\pm}_{\phi}(\alpha)} \prod_{\phi
\neq z, z \pm 1} x_{\phi,\lambda_{\phi}(\alpha)}\\
& = &
\prod_{\phi=z \pm 1} \sum_{\lambda^{\pm}} x_{\phi,\lambda^{\pm}}
\frac{u^{|\lambda^{\pm}|}}{c_{Sp,\phi,q^{m_{\phi}}}(\lambda^{\pm})}  \prod_{\phi
= \bar{\phi} \atop \phi \neq z \pm 1}  \sum_{\lambda}
x_{\phi,\lambda}
\frac{(-(u^{m_{\phi}}))^{|\lambda|}}{c_{GL,z-1,-\sqrt{q^{m_{\phi}}}}(\lambda)}\\
& & \prod_{\{\phi,\bar{\phi}\}, \phi \neq \bar{\phi}} \sum_{\lambda}
x_{\phi,\lambda} x_{\bar{\phi},\lambda}
\frac{u^{2|\lambda|m_{\phi}}}{c_{GL,z-1,q^{m_{\phi}}}(\lambda)} 
\end{eqnarray*}

\end{theorem}

\begin{proof}
	Consider the coefficients of $u^n$ on both sides when $n$ is
even. Their equality follows from Wall's formulas for conjugacy class
sizes for the symplectic groups. We have also made use of the following
elementary fact:

\begin{eqnarray*}
& & \frac{u^{|\lambda_{\phi}|m_{\phi}}}{q^{\frac{m_{\phi}}{2}[2 \sum_{h<i}
hm_h(\lambda_{\phi})m_i(\lambda_{\phi}) + \sum_i
(i-1)m_i(\lambda_{\phi})^2]} \prod_i |U(m_i(\lambda_{\phi}),q^{\frac{m_{\phi}}{2}})|}\\
& = & \frac{(-(u^{m_{\phi}}))^{|\lambda_{\phi}|}} { (-(q^{\frac{m_{\phi}}{2}}))^{[2 \sum_{h<i}
hm_h(\lambda_{\phi})m_i(\lambda_{\phi}) + \sum_i
(i-1)m_i(\lambda_{\phi})^2]} \prod_i
|GL(m_i(\lambda_{\phi}),-(q^{\frac{m_{\phi}}{2}}))|}
\end{eqnarray*}

	which is true because $|U(n,q)| = (-1)^n |GL(n,-q)|$.

	Consider the coefficients of $u^n$ on both sides when $n$ is
odd. The coefficient on the left-hand side is 0. It thus suffices to
show that only even powers of $u$ appear in each term of the product
on the right-hand side of the factorization. This is clear for
polynomials $\phi$ such that $\phi \neq \bar{\phi}$. It is true for
the polynomials $z \pm 1$ because all odd parts in the associated
signed partitions have even multiplicity. Finally, Theorem
$\ref{Count2}$ implies that all polynomials $\phi \neq z \pm 1$ such
that $\phi = \bar{\phi}$ have even degree. So these polynomials
contribute only even powers of $u$ as well.
\end{proof}

\begin{center}
{\bf The Orthogonal Groups}
\end{center}

	This article assumes for simplicity that the characteristic of $F_q$
is not equal to 2. The orthogonal groups can be defined as subgroups
of $GL(n,q)$ preserving a non-degenerate symmetric bilinear form (see
Chapter 1 of Carter $\cite{Ca1}$). For $n=2l+1$ odd, there are two
such forms up to isomorphism, with inner product matrices $A$ and
$\delta A$, where $\delta$ is a non-square in $F_q$ and $A$ is equal
to:

\[ \left( \begin{array}{c c c}
		1 & 0 & 0 \\
		0 & 0_l & I_l \\
		0 & I_l & 0_l
	  \end{array} \right) \]

	Denote the corresponding orthogonal groups by $O^+(2l+1,q)$
and $O^-(2l+1,q)$. This distinction will be useful, even though these
groups are isomorphic. Their common order is:

\[ 2q^{l^2} \prod_{i=1}^l (q^{2i}-1) \]

	For $n=2l$ even, there are again two non-degenerate symmetric
bilinear forms up to isomorphism with inner product matrices:

\[  \left( \begin{array}{c c}
		0_l & I_l \\
		I_l & 0_l
	  \end{array} \right) \] 

\[    \left( \begin{array}{c c c c}
		0_{l-1} & I_{l-1} & 0 & 0 \\
		I_{l-1} & 0_{l-1} & 0 & 0 \\
		0 & 0 & 1 & 0 \\
		0 & 0 & 0 & -\delta
	  \end{array} \right)	\]

	where $\delta$ is a non-square in $F_q$. Denote the
corresponding orthogonal groups by $O^+(2l,q)$ and $O^-(2l,q)$. They
are not isomorphic and have orders:

\[ 2q^{l^2-l} (q^l \mp 1) \prod_{i=1}^{l-1} (q^{2i}-1) \]  

	To describe the conjugacy classes of the finite orthogonal
groups, it is necessary to use the notion of the Witt type of a
non-degenerate quadratic form, as in Chapter 9 of Bourbaki
$\cite{Bou}$. Call a non-degenerate form $N$ null if the vector space
$V$ on which it acts can be written as a direct sum of 2 totally
isotropic subspaces (a totally isotropic space is one on which the
inner product vanishes identically). Define two non-degenerate
quadratic forms $Q'$ and $Q$ to be equivalent if $Q'$ is isomorphic to
the direct sum of $Q$ and a null $N$. The Witt type of $Q$ is the
equivalence class of $Q$ under this equivalence relation. There are 4
Witt types over $F_q$, which Wall denotes by ${\bf 0},{\bf 1}, {\bf
\delta}, {\bf \omega}$, corresponding to the forms $0,x^2,\delta
x^2,x^2 - \delta y^2$ where $\delta$ is a fixed non-square of
$F_q$. These 4 Witt types form a ring, but only the additive structure
is relevant here. The sum of two Witt types with representatives
$Q_1,Q_2$ on $V_1,V_2$ is the equivalence class of $Q_1+Q_2$ on
$V_1+V_2$.

\begin{prop} \label{wit} The four orthogonal groups
$O^+(2n+1,q),O^-(2n+1,q), O^+(2n,q),O^-(2n,q)$ arise from forms $Q$ of
Witt types ${\bf 1},{\bf \delta},{\bf 0},{\bf \omega}$ respectively.
\end{prop}

\begin{proof}
	This follows from the explicit description above of the inner
product matrices which give rise to the various orthogonal groups.
\end{proof}

	Consider the following combinatorial data. To each monic,
non-constant, irreducible polynomial $\phi \neq z \pm 1$ associate a
partition $\lambda_{\phi}$ of some non-negative integer
$|\lambda_{\phi}|$. To $\phi$ equal to $z-1$ or $z+1$ associate an
orthogonal signed partition $\lambda_{\phi}^{\pm}$, by which is meant
a partition of some natural number $|\lambda_{\phi}^{\pm}|$ such that
all even parts have even multiplicity, and all odd $i>0$ have a choice
of sign. For $\phi= z-1$ or $\phi = z+1$ and odd $i>0$, we denote by
$\Theta_i (\lambda_{\phi}^{\pm})$ the Witt type of the orthogonal
group on a vector space of dimension $m_i(\lambda_{\phi}^{\pm})$ and
sign the choice of sign for $i$.

\begin{center}
Example of an Orthogonal Signed Partition
\end{center}
	
\[ \begin{array}{c c c c c}
	&  . & . & . & .   \\
	&  . & . & . & .   \\
	-  & . & . & . &   \\
	&  . & . &  &   \\
	&  . & . &  &    \\
	+  & . &  & &   \\
	 & . &  & &   
	  \end{array}  \]

	Here the $-$ corresponds to the part of size 3 and the $+$
corresponds to the parts of size 1. 

	Theorem \ref{almostWall}, though not explicitly stated by Wall
\cite{Wal}, is implicit in the discussion on pages 38-40 of his article. His
statements seem different because he fixes an
$\alpha \in GL(n,q)$ and asks which orthogonal groups contain a conjugate of
$\alpha$. Here we want to fix the group and parameterize its conjugacy
classes. The polynomial $\bar{\phi}$ is defined as in the symplectic case.

\begin{theorem} \label{almostWall} The data $\lambda^{\pm}_{z-1},
\lambda^{\pm}_{z+1}, \lambda_{\phi}$ represents a conjugacy class of
some orthogonal group if:

\begin{enumerate}
\item $|\lambda_{z}|=0$
\item $\lambda_{\phi}=\lambda_{\bar{\phi}}$
\item $\sum_{\phi=z \pm 1} |\lambda_{\phi}^{\pm}| + \sum_{\phi \neq z \pm
1} |\lambda_{\phi}| m_{\phi}=n$
\end{enumerate}
	
	In this case, the data represents the conjugacy class of
exactly 1 orthogonal group $O(n,q)$, with sign determined by the
condition that the group arises as the stabilizer of a form of Witt
type:

\[ \sum_{\phi=z \pm 1} \sum_{i \ odd} \Theta_i(\lambda_{\phi}^{\pm}) +
\sum_{\phi \neq z \pm 1} \sum_{i \geq 1} i m_i(\lambda_{\phi}) {\bf \omega} \]

	The conjugacy class has size:

\[ \frac{|O(n,q)|}{\prod_{\phi} B(\phi)} \]

	where

\begin{eqnarray*}
B(\phi) & = & q^{[\sum_{h<i}
hm_h(\lambda_{\phi}^{\pm})m_i(\lambda_{\phi}^{\pm}) + \frac{1}{2} \sum_i
(i-1)m_i(\lambda_{\phi}^{\pm})^2]} \prod_i A(\phi^{\pm,i}) \ if \ \phi=z
\pm 1\\
B(\phi) & = & q^{m_{\phi}[\sum_{h<i}
hm_h(\lambda_{\phi})m_i(\lambda_{\phi}) + \frac{1}{2} \sum_i
(i-1)m_i(\lambda_{\phi})^2]} \prod_i A(\phi^i) \ if \ \phi \neq z \pm 1
\end{eqnarray*}
  
	and

\begin{eqnarray*}
A(\phi^{\pm,i}) & = & |O(m_i(\lambda_{\phi}^{\pm}),q)| \ if
\ i=1 \ (mod \ 2)\\
& = & q^{-\frac{m_i(\lambda_{\phi}^{\pm})}{2}}
|Sp(m_i(\lambda_{\phi}^{\pm}),q)| \ if \ i=0 \ (mod \ 2)\\
A(\phi^i) & = & |U(m_i(\lambda_{\phi}),q^{\frac{m_{\lambda}}{2}})| \ if \
\phi=\bar{\phi}\\
& = & |GL(m_i(\lambda_{\phi}),q^{m_{\lambda}})|^{\frac{1}{2}} \ if \ \phi \neq
\bar{\phi}.
\end{eqnarray*}

	Here $O(m_i(\lambda_{\phi}),q)$ is the orthogonal group with
the same sign as the sign associated to the parts of size $i$.
\end{theorem}
 
	In the case $\phi= z \pm 1$, $B(\phi)$ will also be denoted by
$c_{O,\phi,q^{m_{\phi}}}(\lambda^{\pm})$.

	As an example of this formula, consider the set of orthogonal
symmetries in $O^+(n,q)$ where $n$ is odd (considerations for the
other orthogonal groups are analogous). An orthogonal symmetry is a
determinant $-1$ orthogonal map with an $n-1$ dimensional fixed
space. Orthogonal symmetries are important because they generate the
orthogonal group containing them (page 129 of Artin
$\cite{Artin}$). It is worth remarking that orthogonal transvections
exist only in the characteristic 2 case, which is excluded here.

	Fulman \cite{fulthesis} counts the orthogonal symmetries in $O^+(n,q)$
directly, showing that there are $q^{n-1}$ of them. The orthogonal symmetries in
$O^+(n,q)$ fall into two conjugacy classes. These may be described in terms of
Wall's combinatorial data. One conjugacy class corresponds to the data
$\lambda_{z-1}=(+ 1^{n-1}), \lambda_{z+1} = (+1)$. The other conjugacy
class corresponds to the data $\lambda_{z-1} = (-1^{n-1}),
\lambda_{z+1} = (-1)$. This follows from Theorem $\ref{almostWall}$
and the fact that the dimension of the fixed
space of $\alpha$ is the number of parts of $\lambda^{\pm}_{z-1}
(\alpha)$. So the class size formulas of Wall imply that the total
number of orthogonal symmetries in $O^+(n,q)$ for $n$ odd is:

\[ \frac{|O^+(n,q)|}{q^{\frac{1}{2}(n-2)} |O^+(n-1,q)|
q^{\frac{1}{2}} |O^+(1,q)|} + \frac{|O^+(n,q)|}{q^{\frac{1}{2}(n-2)} |O^-(n-1,q)|
q^{\frac{1}{2}} |O^-(1,q)|} = q^{n-1} \]

	In analogy with the general linear, unitary, and symplectic groups, define a
cycle index for the orthogonal groups by:

\begin{eqnarray*}
\sum_{n=1}^{\infty} [\frac{u^n}{|O^+(n,q)|} \sum_{\alpha \in O^+(n,q)}
\prod_{\phi=z \pm 1} x_{\phi,\lambda_{\phi}^{\pm}(\alpha)}
\prod_{\phi \neq z, z \pm 1} x_{\phi,\lambda_{\phi}(\alpha)}]\\
+ [\frac{u^n}{|O^-(n,q)|} \sum_{\alpha \in O^-(n,q)} \prod_{\phi=z \pm 1}
x_{\phi,\lambda_{\phi}^{\pm}(\alpha)} \prod_{\phi \neq z, z \pm 1}
x_{\phi,\lambda_{\phi}(\alpha)}]
\end{eqnarray*}

	Theorem \ref{cyort} shows that the cycle index of the finite orthogonal groups
factors.

\begin{theorem} \label{cyort}

\begin{eqnarray*}
& & \sum_{n=1}^{\infty} [\frac{u^n}{|O^+(n,q)|} \sum_{\alpha \in O^+(n,q)}
\prod_{\phi=z \pm 1} x_{\phi,\lambda_{\phi}^{\pm}(\alpha)}
\prod_{\phi \neq z, z \pm 1} x_{\phi,\lambda_{\phi}(\alpha)}]\\
& & + [\frac{u^n}{|O^-(n,q)|} \sum_{\alpha \in O^-(n,q)} \prod_{\phi=z \pm 1}
x_{\phi,\lambda_{\phi}^{\pm}(\alpha)} \prod_{\phi \neq z, z \pm 1}
x_{\phi,\lambda_{\phi}(\alpha)}]\\
& = & \prod_{\phi=z \pm 1}
\sum_{\lambda^{\pm}} x_{\phi,\lambda^{\pm}}
\frac{u^{|\lambda^{\pm}|}}{c_{O,\phi,q^{m_{\phi}}}(\lambda^{\pm})}
 \prod_{\phi
= \bar{\phi} \atop \phi \neq z \pm 1}  \sum_{\lambda}
x_{\phi,\lambda}
\frac{(-(u^{m_{\phi}}))^{|\lambda|}}{c_{GL,z-1,-\sqrt{q^{m_{\phi}}}}(\lambda)}\\
& & \prod_{\{\phi,\bar{\phi}\}, \phi \neq \bar{\phi}} \sum_{\lambda}
x_{\phi,\lambda} x_{\bar{\phi},\lambda}
\frac{u^{2|\lambda|m_{\phi}}}{c_{GL,z-1,q^{m_{\phi}}}(\lambda)}
\end{eqnarray*}

\end{theorem}

\begin{proof}
	The first equation follows from Theorem $\ref{almostWall}$
since every term in the product on the right hand side corresponds to
a conjugacy class in exactly one of the orthogonal groups, and the
class sizes check.
\end{proof}

\section{The $q \rightarrow \infty$ Limit of the Cycle Indices}
\label{BIGQ}

	At the end of his paper, Stong $\cite{St1}$ stated that if one sets
$x_{\phi,\lambda}=(x_{m_{\phi}})^{|\lambda|}$ in the cycle index of the general
linear groups and lets $q \rightarrow \infty$, then one obtains the cycle index of the symmetric
groups. Indeed, Lemma \ref{Sto} from Section $\ref{GLCYC}$ shows that for fixed
$q$ the cycle index for $GL$ becomes:

\[ \prod_{m=1}^{\infty} \prod_{r=1}^{\infty}
(\frac{1}{1-(\frac{u}{q^r})^m x_m})^{I_{m,q}} \]

	where $I_{m,q}$ is the number of monic, degree $m$, irreducible $\phi \neq
z$ defined over a field of $q$ elements. Letting $q \rightarrow \infty$ and using
the formula for
$I_{m,q}$ (Lemma $\ref{irred}$) gives:

\[ lim_{q \rightarrow \infty} \prod_{m=1}^{\infty}
\prod_{r=1}^{\infty} (\frac{1}{1-(\frac{u}{q^r})^m x})^{I_{m,q}} =
\prod_{m=1}^{\infty} e^{\frac{x_mu^m}{m}} \]

	This fact may also be stated as follows. Fix $n$. Then the $q
\rightarrow \infty$ limit of the probability that the characteristic
polynomial of a uniformly chosen element of $GL(n,q)$ factors into $a_m$
irreducible polynomials of degree $m$ is equal to the chance that a randomly
chosen element of $S_n$ has $a_m$ $m$-cycles.

	We now use algebraic groups to give a conceptual statement and
proof of Stong's observation which generalizes to other groups. Dick
Gross suggested that such an interpretation should exist and was kind
enough to explain the basics of algebraic groups.

	The necessary background about maximal tori in algebraic
groups can all be found in Chapter 3 of Carter $\cite{Ca2}$. Let us
review some of these facts. Take $G$ to be a connected, reductive
algebraic group over $\bar{F_q}$ which is Chevalley (this means that
the Frobenius map giving rise to $G^F$ is $x \rightarrow x^q$ where
$q$ is a prime power) such that $G'$ is simply connected. Recall that
a torus of $G$ is a subgroup of $G$ which is isomorphic to a product
of copies of the multiplicative group of $\bar{F_q}$, and that a
maximal torus of $G$ is a maximal such subgroup. Any two maximal tori
of $G$ are conjugate in $G$.

	A maximal torus of $G^F$ is defined to be a group $T^F$ where
$T$ is a maximal torus of $G$. As is discussed on pages 32-33 of
Carter, the Lang-Steinberg theorem of algebraic groups implies that
$G^F$ has a maximal torus $T^F$ which is diagonalizable over $F_q$
(such a maximal torus is called maximally split). While it is true
that all maximally split maximal tori $T^F$ are conjugate in $G^F$,
the maximal tori of $G^F$ may fall into many conjugacy classes.

	Proposition 3.3.3 of Carter says that under these conditions,
there is a bijection $\Phi$ between $G^F$ conjugacy classes of
$F$-stable maximal tori in $G$ and conjugacy classes of the Weyl group
$W$. We recall the definition of $\Phi$ (the proof that it is a
bijection is harder). Let $T_0$ be a fixed $F$-stable maximally split
maximal torus of $G$ (this exists by the Lang-Steinberg
theorem). Since all maximal tori in $G$ are conjugate to $T_0$, one
can write $T=^gT_0$ for some $g \in G$ (the symbol $^g$ denotes
conjugation by $g$). Clearly $g^{-1}F(g) \in N(T_0)$. Since
$W=N(T_0)/T$, this associates to $T$ an element of $W$, which turns
out to be well defined up to conjugacy in $W$.

	Next, one can define a map $\omega$ (similar to that in Lehrer
$\cite{Leh}$) from $G^F$ to conjugacy classes of $W$. Given $\alpha
\in G^F$, let $\alpha_s$ be the semi-simple part of $\alpha$. Theorem
3.5.6 of Carter says that $G'$ simply connected implies that
$C_G(\alpha_s)$ is connected. Take $T$ to be an $F$-stable maximal
torus in $C_G(\alpha_s)$ such that $T^F$ is maximally split. By what
has been said before, all such $T$ are conjugate in $C_G(\alpha_s)^F$
and hence in $G^F$. Define $\omega(\alpha)=\Phi(T)$.

	In the case of $GL(n,q)$, the map $\omega$ sends an element
$\alpha$ whose characteristic polynomial factors into $a_m$
irreducible polynomials of degree $m$ to the conjugacy class of $S_n$
corresponding to permutations with $a_m$ $m$-cycles.

	Proposition 3.6.6 of Carter implies that for $q$ sufficiently
large, all maximal tori $T^F$ of $G^F$ lie in exactly 1 maximal torus
of $G$ (this condition is called non-degeneracy). For such $q$ one can
then define a a bijection $\Phi'$ between $G^F$ conjugacy classes of
maximal tori $T^F$ of $G^F$ and conjugacy classes of $W$ by
$\Phi'(T^F)=\Phi(T)$, where $T$ is the unique maximal torus of $G$
containing $T^F$.

	We now prove the following theorem.

\begin{theorem} \label{tori} Let $G$ be a connected, reductive
Chevalley group which is defined over $F_q$, such that $G'$ is simply
connected. Suppose that as $q \rightarrow \infty$, the chance that an
element of $G^F$ is regular, semi-simple approaches 1. Then for all
conjugacy classes $c$ in $W$,

\[ \lim_{q \rightarrow \infty} P_{G^F}(\omega(\alpha) \in  c)=P_W(w \in c) \]  

where both probabilities are with respect to the uniform distribution.
\end{theorem}

\begin{proof}
	Take $q$ large enough that all maximal tori of $G^F$ are
non-degenerate, so that the construction of $\Phi'$ works. From page
29 of Carter $\cite{Ca2}$, a regular semi-simple element $\alpha$ of
$G$ lies in a unique maximal torus, which implies by non-degeneracy
that $\alpha$ lies in a unique $T^F$. This also implies that
$\omega(\alpha)=\Phi'(T^F)$. Therefore:

\begin{eqnarray*}
\lim_{q \rightarrow \infty} P_{G^F}(\omega(\alpha)=c) & = & \lim_{q
\rightarrow \infty} \frac{|\{\alpha \ regular \ semisimple :
\omega(\alpha)=c\}|}{|G^F|}\\
& = & \lim_{q \rightarrow \infty} \sum_{T^F:\Phi'(T^F)=c}  \frac{|\{\alpha \ regular \ semisimple, \alpha \in T^F \}|}{|G^F|}\\
& = &  \lim_{q \rightarrow \infty} \frac{|\{\alpha \ regular \ semisimple, \alpha \in T\}|}{|N_{G^F}(T^F)|}
\end{eqnarray*}

	Proposition 3.3.6 and Corollary 3.6.5 of Carter give that
$N_{G^F}(T^F)/T^F$ is isomorphic to $C_W(w)$, so that:

	\[ |N_{G^F}(T^F)|/|T^F| = \frac{|W|}{|c|} \]

	Therefore,

\[ \lim_{q \rightarrow \infty} P_{G^F}(\omega(\alpha)=c) =
\frac{|c|}{|W|} \lim_{q \rightarrow \infty} \frac{|\{\alpha \
regular \ semisimple, \alpha \in T^F\}|}{|T^F|} \]

	Summing over the finitely many conjugacy classes $c$ on both sides of
this equation gives 1. Thus,

\[ \lim_{q \rightarrow \infty} \frac{|\{\alpha \
regular \ semisimple, \alpha \in T^F\}|}{|T^F|} = 1 \]

	for all $T^F$, which proves the theorem.
\end{proof}

	The heuristics at the end of Section $\ref{REG}$ suggest
that as $q \rightarrow \infty$, the chance that an element is regular
semi-simple approaches 1. This should be directly checkable for the
classical groups using the cycle indices (although the
unitary group is not Chevalley so Theorem $\ref{tori}$ does not
apply). Some recent work of Lehrer $\cite{Le2}$ using $l$-adic
cohomology gives involved expressions for the chance that an element
of a finite group of Lie type is regular semisimple. It may be
possible to read the $q \rightarrow \infty$ limit off of his results.

\section{Application: The Characteristic Polynomial, Number of Jordan
Blocks, and Average Order of an Element in a Finite Classical Group}
\label{APPLICATIONS}

	This section applies the cycle indices of Section \ref{OTHER} to carry over work
of Stong
\cite{St1},\cite{St2} on the general linear groups to the unitary,
symplectic, and orthogonal groups. In particular,
we will study the characteristic polynomial, the number of Jordan blocks, and the
average order of an element in a finite classical group. Other authors, notably
Hansen \cite{Han} and Goh and Schmutz \cite{Go2}, have obtained interesting
results about random matrices. The results of this section are an indication that
the analogs of their theorems carry over to the other classical groups.

\begin{center}
{\bf Application 1: The Characteristic Polynomial}
\end{center}

	This application will make use of the following theorem of Steinberg,
which counts the number of unipotent elements (i.e. all eigenvalues equal to one)
$\alpha$ in a finite classical group. Clearly $\alpha$ is unipotent if and only
if its characteristic polynomial is a power of $z-1$. Steinberg's count is
normally proven using the Steinberg character, as on page 156 of Humphreys
$\cite{Hum}$. Fulman \cite{fulprob} gives a probabilistic proof, at least
for the general linear and unitary cases.

\begin{theorem} \label{Steinberg} The number of unipotent elements in a
finite group of Lie type $G^F$ is the square of the order of a $p$-Sylow of
$G^F$, where $p$ is the prime used in the construction of $G^F$ (in the
case of the classical groups, $p$ is the characteristic of $F_q$).
\end{theorem}

{\bf The General Linear Groups} Reiner \cite{Reiner} and Gerstenhaber
\cite{Ger} proved the following theorem counting elements of $GL(n,q)$
wight a given characteristic polynomial. Presumably Stong \cite{St1} knew
how to do this using the cycle index, but we include a proof as one did not
appear there.

	Recall that if $f(u)$ is a polynomial in $u^n$, then the notation $[u^n] f(u)$
means the coefficient of $u^n$ in $f(u)$.

\begin{theorem} Let $\phi$ be a monic polynomial of degree $n$
which factors into irreducibles as $\phi= \prod_{i=1}^r
\phi_i^{j_i}$. Then the number of elements of $GL(n,q)$ with
characteristic polynomial $\phi$ is:

\[ |GL(n,q)| \prod_{i=1}^r
\frac{q^{m_{\phi}j_i(j_i-1)}}{|GL(j_i,q^{m_{\phi}})|} \]
\end{theorem}

\begin{proof}
	In the cycle index for the general linear groups, set $x_{\phi,\lambda}=1$ if
$\phi=\phi_i$ and $|\lambda|=j_i$. Otherwise, set $x_{\phi,\lambda}=0$. Taking
coefficients of $u^n$ shows that the number of elements of $GL(n,q)$ with
characteristic polynomial $\phi$ is:

\begin{eqnarray*}
& & |GL(n,q)| [u^n] \prod_{i=1}^r \sum_{\lambda: |\lambda|=j_i}
\frac{u^{j_i(m_{\phi_i})}}{c_{GL,\phi,q^{m_{\phi}}}}\\
& = & |GL(n,q)| \prod_{i=1}^r \sum_{\lambda: |\lambda|=j_i}
\frac{1}{c_{GL,\phi,q^{m_{\phi}}}}\\
& = & |GL(n,q)| \prod_{i=1}^r
\frac{q^{m_{\phi}j_i(j_i-1)}}{|GL(j_i,q^{m_{\phi}})|}
\end{eqnarray*}

	where the last equality uses the cycle index and Steinberg's count of
unipotents (Theorem \ref{Steinberg}) which says that there are $q^{n(n-1)}$
unipotent elements in $GL(n,q)$.
\end{proof}

{\bf The Unitary Groups} Arguing as in the general linear case, but using the
cycle index for the unitary groups leads to the following count of elements of
$U(n,q)$ given a given characteristic polynomial.

\begin{theorem} \label{Uchar} Let $\phi$ be a monic polynomial of
degree $n$ which factors into irreducibles as $\phi=\prod_i
\phi_i^{j_i} \prod_{i'} [\phi_{i'} \tilde{\phi_{i'}}]^{j_{i'}}$, where
$\phi_i=\tilde{\phi_i}$ and $\phi_{i'} \neq \tilde{\phi_{i'}}$. Then
the number of elements of $U(n,q)$ with characteristic polynomial
$\phi$ is:

\[ |U(n,q)| \prod_i \frac{q^{m_{\phi_i} j_i(j_i-1)}}
{|U(j_i,q^{m_{\phi_i}})|} \prod_{i'}
\frac{q^{2m_{\phi_{i'}}j_{i'}(j_{i'}-1)}}{|GL(j_{i'},
q^{2m_{\phi_{i'}}})|} \]

\end{theorem}

{\bf The Symplectic Groups} The statement for the symplectic groups, proved by
the same technique as for the general linear and unitary groups, is as follows.

\begin{theorem} \label{SpChar} Let $\phi$ be a polynomial of degree 2n
which factors into irreducibles as $(z-1)^{2a}(z+1)^{2b} \prod_{i}
\phi_i^{j_i} \prod_{i'} [\phi_{i'} \bar{\phi_{i'}}]^{j_{i'}}$ where
$\phi_{i'} \neq \bar{\phi_{i'}}$. Then the number of elements of $Sp(2n,q)$
with characteristic polynomial $\phi$ is:

\[ |Sp(2n,q)| \frac{q^{2a^2}}{|Sp(2a,q)|} \frac{q^{2b^2}}{|Sp(2b,q)|}
\prod_{i} \frac{q^{\frac{m_{\phi_i}j_i(j_i-1)}{2}}}{|U(j_i,q^{\frac{m_
{\phi_i}}{2}})|} \prod_{i'} \frac{q^{m_{\phi_{i'}}j_i(j_i-1)}}
{|GL(j_i,q^{m_{\phi_{i'}}})|} \]

\end{theorem}

{\bf The Orthogonal Groups} The corresponding statement for the orthogonal groups
is somewhat more complicated, because the cycle index involved averaging over
$O^+(n,q)$ and $O^-(n,q)$.

\begin{theorem} \label{Ocharpoly} Let $\phi$ be a polynomial of degree $n$
which factors into irreducibles as $(z-1)^a (z+1)^b \prod_i \phi_i^{j_i}
\prod_{i'} [\phi_{i'} \bar{\phi_{i'}}]^{j_{i'}}$ where $\phi_{i'} \neq
\bar{\phi_{i'}}$. Then $\frac{1}{2}$ of the sum of the proportion of
elements in $O^+(n,q)$ and $O^-(n,q)$ with characteristic polynomial $\phi$
is:

\[ \frac{F(a) F(b)}{2} \prod_{i}
\frac{q^{\frac{m_{\phi_i}j_i(j_i-1)}{2}}}{|U(j_i,q^{\frac{m_
{\phi_i}}{2}})|} \prod_{i'} \frac{q^{m_{\phi_{i'}}j_i(j_i-1)}}
{|GL(j_i,q^{m_{\phi_{i'}}})|} \]

	where:

\begin{eqnarray*}
F(n) & = & \frac{q^{\frac{n^2}{2}}}{|Sp(n,q)|} \ if \ n = 0 \ (mod \ 2)\\
 & = & \frac{q^{\frac{(n-1)^2}{2}}}{|Sp(n-1,q)|} \ if \ n = 1 \ (mod \ 2)
\end{eqnarray*}
\end{theorem}

\begin{center}
{\bf Application 2: Number of Jordan Blocks}
\end{center}

	For $\alpha \in Mat(n,q)$, let $X_n(\alpha)$ be the number of
irreducible polynomials counted with multiplicity occurring in the
rational canonical form of $\alpha$. (This is not quite the number of Jordan
blocks, but is in the limit that $q \rightarrow \infty$). This application uses
the cycle index to study $X_n(\alpha)$ where $\alpha$ is an element of a finite
classical group. In fact only the mean of $X_n$ will be computed, but the more
analytically inclined should be able to compute the variance and prove asymptotic
normality using the cycle index.

{\bf General Linear Groups} The work has been done already in this case. Stong
$\cite{St1}$ proves that the random variable $X_n$ has mean and variance
$log(n)+O(1)$. Goh and Schmutz $\cite{Go2}$ prove that $X_n$ is asymptotically
normal.

	We give a somewhat simpler computation of the mean, both to illustrate the
elegance of the cycle index approach, and because the same technique works for the
other classical groups.

\begin{lemma} \label{jord}
 \[ \sum_{n=0}^{\infty} (1-u)u^n \sum_{\alpha \in GL(n,q)}
x^{X_n(\alpha)} = \prod_{m=1}^{\infty} \prod_{i=1}^{\infty}
(\frac{1-(\frac{u}{q^i})^m}{1-(\frac{u}{q^i})^m x})^{I_{m,q}} \]
\end{lemma}

\begin{proof}
	Set $x_{\phi,\lambda}=x^{|\lambda|}$ for all polynomials
$\phi$. The result follows from Lemma \ref{Sto} in Section \ref{GLCYC}.
\end{proof}

	The mean of $X_n$ is now easily computed.

\begin{theorem} $EX_n = log(n) + O(1)$, where the expectation is taken
over the group $GL(n,q)$ with $q$ fixed.
\end{theorem}

\begin{proof}
	Differentiating both sides of the generating function of Lemma
$\ref{jord}$ with respect to $x$ and then setting $x=1$ gives:

\begin{eqnarray*}
EX_n & = & [u^n] \frac{1}{1-u} \sum_{m=1}^{\infty} I_{m,q} \sum_{i=1}^{\infty} \sum_{l=1}^{\infty} \frac{u^{ml}}{q^{iml}}\\
& = & \sum_{r=1}^n (\sum_{m|r} I_{m,q}) (\sum_{i=1}^{\infty} \frac{1}{q^{ri}})
\end{eqnarray*}

	It is well known that $I_{m,q}=\frac{q^m}{m} +
O(q^{\frac{m}{2}})$, from which it follows that $\sum_{m|r}
I_{m,q}=\frac{q^r}{r} + O(q^{\frac{r}{2}})$. Therefore:

\begin{eqnarray*}
EX_n & = & \sum_{r=1}^n \sum_{i=1}^{\infty} \frac{1}{q^{ri}} (\frac{q^r}{r} + O(q^{\frac{r}{2}}))\\
& = & \sum_{r=1}^n \frac{1}{r} + O(q^{-\frac{r}{2}})\\
& = & (\sum_{r=1}^n \frac{1}{r}) + O(1)\\
& = & log(n) + O(1)
\end{eqnarray*}

\end{proof}

{\bf The Unitary Groups} The same technique works for the unitary groups.

\begin{lemma} \label{unipgenfun}

\[ \sum_{n=0}^{\infty} (1-u)u^n \sum_{\alpha \in U(n,q)}
x^{X_n(\alpha)} = \prod_{m=1 \atop m \ odd}^{\infty}
\prod_{i=1}^{\infty}
(\frac{1+(-1)^i(\frac{u^m}{q^{im}})}{1+(-1)^i(\frac{u^m}{q^{im}})x})^{\tilde{I}_{m,q^2}}
\prod_{m=1}^{\infty} \prod_{i=1}^{\infty}
(\frac{1-(\frac{u^{2m}}{q^{2im}})}{1-(\frac{u^{2m}}{q^{2im}})x^2})^{\frac{I_{m,q^2}-\tilde{I}_{m,q^2}}{2}}
\]

\end{lemma}

\begin{proof}
	Proceed as in Lemma \ref{jord}.
\end{proof}

\begin{theorem} $EX_n = \frac{3}{2}log(n) + O(1)$, where the expectation
is taken over the group $U(n,q)$ with $q$ fixed.  \end{theorem}

\begin{proof}
	Differentiating both sides of the generating function of Lemma
$\ref{unipgenfun}$ with respect to $x$ and setting $x=1$ gives:

\begin{eqnarray*}
EX_n & = & [u^n] \frac{1}{1-u} [\sum_{m=1 \atop m \ odd}^{\infty}
\tilde{I}_{m,q^2} \sum_{i=1}^{\infty} \sum_{l=1}^{\infty} (-1)^{(i+1)(l+1)}
\frac{u^{ml}}{q^{iml}}] + [\sum_{m=1}^{\infty}
(I_{m,q^2}-\tilde{I}_{m,q^2}) \sum_{i=1}^{\infty} \sum_{l=1}^{\infty}
\frac{u^{2ml}}{q^{2iml}}]\\
& = & [\sum_{r=1}^n \sum_{m|r \atop m \ odd} \tilde{I}_{m,q^2}
\sum_{i=1}^{\infty} (-1)^{(i+1)(\frac{r}{m}+1)} \frac{1}{q^{ri}}] +
[\sum_{r=1 \atop r \ even}^{n} \sum_{m | \frac{r}{2}} (I_{m,q^2}-\tilde{I}_{m,q^2}) \sum_{i=1}^{\infty}
\frac{1}{q^{ri}} ] 
\end{eqnarray*}

	Recall from Theorem $\ref{Count1}$ that $\tilde{I}_{m,q^2} =
\frac{q^m}{m} + O(q^{\frac{m}{2}})$ for $m$ odd. We also know that
$I_{m,q^2}=\frac{q^{2m}}{m} + O(q^m)$. Thus the dominant contribution from
the first bracketed term comes from $m=r,i=1,l=1$, and the dominant
contribution from the second bracketed term comes from
$m=\frac{r}{2},i=1,l=1$. Therefore,

\begin{eqnarray*}
EX_n & = & [\sum_{r=1 \atop r \ odd}^n \frac{1}{r}+O(1)] + [\sum_{r=1 \atop
r \ even}^n \frac{1}{\frac{r}{2}}+O(1)]\\
& = & \sum_{r=1}^n \frac{1}{r} + \frac{1}{2} \sum_{r=1}^{\frac{n}{2}} \frac{1}{r} + O(1)\\
& = & \frac{3}{2} log(n) + O(1)
\end{eqnarray*}

\end{proof}

{\bf The Symplectic Groups} As with the general linear and unitary groups, let
$X_{2n}(\alpha)$ be the number of irreducible polynomials counted with
multiplicity in the Jordan canonical form of $\alpha \in
Sp(2n,q)$. Lemma $\ref{Symge}$ gives a generating function for
$X_{2n}(\alpha)$.

\begin{lemma} \label{Symge}

\begin{eqnarray*}
\sum_{n=0}^{\infty} (1-u)u^{2n} \frac{1}{|Sp(2n,q)|} \sum_{\alpha \in Sp(2n,q)}
x^{X_{2n}(\alpha)} & = & \prod_{i=1}^{\infty}
(\frac{1-\frac{u^2}{q^{2i-1}}} {1-\frac{(ux)^2}{q^{2i-1}}})^2 \prod_{m=1
\atop m \ even}^{\infty} \prod_{i=1}^{\infty}
(\frac{1+(-1)^i(\frac{u^m}{q^{\frac{im}{2}}})}
{1+(-1)^i(\frac{(ux)^m}{q^{\frac{im}{2}}})})^{\bar{I}_{m,q}}\\
&  & \prod_{m=1}^{\infty}
\prod_{i=1}^{\infty} (\frac{1-\frac{u^{2m}}{q^{im}}}
{1-\frac{(ux)^{2m}}{q^{im}}})^{\frac{I_{m,q}-\bar{I}_{m,q}}{2}}
\end{eqnarray*}

\end{lemma}

\begin{proof}
	This follows immediately from the cycle index for the
symplectic groups.
\end{proof}

	With Lemma $\ref{Symge}$, the mean of $X_{2n}$ over $Sp(2n,q)$
is straightforward to compute.

\begin{theorem} \label{JordanSp} $EX_{2n} = \frac{3}{2}logn + O(1)$
\end{theorem}

\begin{proof}
	Differentiating both sides of the generating function of Lemma
$\ref{Symge}$ with respect to $x$ and then setting $x=1$ gives:

\begin{eqnarray*}
EX_{2n} & = & \sum_{r=1}^n [u^{2r}] ((2 \sum_{i=1}^{\infty}
\frac{\frac{2u}{q^{2i-1}}} {1-\frac{u^i}{q^{2i-1}}}) + (\sum_{m=1 \atop m \
even}^{\infty} \bar{I}_{m,q} \sum_{i=1}^{\infty} \sum_{l=1}^{\infty}
(-1)^{(i+1)(l+1)}\frac{u^{ml}}{q^{\frac{iml}{2}}})\\
& = & \ \ \ \ \  + (\sum_{m=1}^{\infty}
(I_{m,q}-\bar{I}_{m,q}) \sum_{i=1}^{\infty} \sum_{l=1}^{\infty}
\frac{u^{2ml}}{q^{iml}}))
\end{eqnarray*}

	The term coming from the first sum is $O(1)$. By Theorem
$\ref{Count2}$, $\bar{I}_{m,q}=\frac{q^{\frac{m}{2}}}{m} +
O(q^{\frac{m}{4}})$ for $m$ even. We also know that
$I_{m,q}=\frac{q^m}{m} + O(q^{\frac{m}{2}})$. The dominant
contribution to $[u^{2r}]$ from the second and third sums comes from
$m=2r,i=1,l=1$ and $m=r,i=1,l=1$ respectively. Thus,
 
\begin{eqnarray*}
EX_{2n} & = & [\sum_{r=1}^{n} \frac{1}{2r}] + [\sum_{r=1}^{n} \frac{1}{r}]
+ O(1)\\
& = & \frac{3}{2}logn + O(1)
\end{eqnarray*}

\end{proof}

{\bf The Orthogonal Groups} As with the other classical groups, let $X_n(\alpha)$
be the number of irreducible polynomials counted with multiplicity in the
Jordan canonical form of $\alpha \in O(n,q)$.

\begin{lemma} \label{genjoror}

\begin{eqnarray*}
\sum_{n=0}^{\infty} (\frac{1-u}{1+u}) u^n [\sum_{\alpha \in O^+(n,q)}
\frac{x^{X_n(\alpha)}}{|O^+(n,q)|} + \sum_{\alpha \in O^-(n,q)}
\frac{x^{X_n(\alpha)}}{|O^-(n,q)|} ] & = &  (\frac{1+ux}{1+u} \prod_{i=1}^{\infty}
\frac{1-\frac{u^2}{q^{2i-1}}} {1-\frac{(ux)^2}{q^{2i-1}}})^2\\
&  & \prod_{m=1
\atop m \ even}^{\infty} \prod_{i=1}^{\infty}
(\frac{1+(-1)^i\frac{u^m}{q^{\frac{im}{2}}}}
{1+(-1)^i\frac{(ux)^m}{q^{\frac{im}{2}}}})^{\bar{I}_{m,q}}\\
&  & \prod_{m=1}^{\infty}
\prod_{i=1}^{\infty} (\frac{1-\frac{u^{2m}}{q^{im}}}
{1-\frac{(ux)^{2m}}{q^{im}}})^{\frac{I_{m,q}-\bar{I}_{m,q}}{2}}
\end{eqnarray*}

\end{lemma}   

\begin{proof}
	This follows from the cycle index for the orthogonal groups.
\end{proof}

	Let $EX_n$ be $\frac{1}{2}$ of the sum of the averages of $X_n$
over $O^+(n,q)$ and $O^-(n,q)$.

\begin{theorem} $EX_n = \frac{3}{2}log(n) + O(1)$
\end{theorem}

\begin{proof}
	To compute $EX_n$ differentiate both sides and set $x=1$ in:

\[ [u^n] \frac{1+u}{2(1-u)} (\frac{1+ux}{1+u} \prod_{i=1}^{\infty}
\frac{1-\frac{u^2}{q^{2i-1}}} {1-\frac{(ux)^2}{q^{2i-1}}})^2 \prod_{m=1
\atop m \ even}^{\infty} \prod_{i=1}^{\infty}
(\frac{1+(-1)^i\frac{u^m}{q^{\frac{im}{2}}}}
{1+(-1)^i\frac{(ux)^m}{q^{\frac{im}{2}}}})^{\bar{I}_{m,q}} \prod_{m=1}^{\infty}
\prod_{i=1}^{\infty} (\frac{1-\frac{u^{2m}}{q^{im}}}
{1-\frac{(ux)^{2m}}{q^{im}}})^{\frac{I_{m,q}-\bar{I}_{m,q}}{2}} \]

	The first product contributes $O(1)$, so that:

\[ EX_n = \sum_{r=1}^n \frac{[u^r]+[u^{r-1}]}{2} ((\sum_{m=1 \atop m \
even}^{\infty} \bar{I}_{m,q} \sum_{i=1}^{\infty} \sum_{l=1}^{\infty}
(-1)^{(i+1)(l+1)} \frac{u^{ml}}{q^{\frac{iml}{2}}}) + (\sum_{m=1}^{\infty}
(I_{m,q}-\bar{I}_{m,q}) \sum_{i=1}^{\infty} \sum_{l=1}^{\infty}
\frac{u^{2ml}}{q^{iml}})) \]
	
	This proves the theorem by the same logic as in Theorem
$\ref{JordanSp}$.
\end{proof}

\begin{center}
{\bf Application 3: The Average Order}
\end{center}

	Let $v_{GL,n}$ be the average over $GL(n,q)$ of the order of
an element of $GL(n,q)$. Stong
$\cite{St2}$ shows that for fixed $q$ and growing $n$, $log(v_{GL,n})
= nlog(q) - log(n) + o(log(n))$. The lower bound in the equality was not
difficult; the upper bound was hard.

	In what follows we compute lower bounds for the other classical groups.
Presumably these bounds give the correct answer, as in the general linear case.
The involutions $\tilde{}$ and $\bar{}$ on polynomials are defined as in Section
\ref{OTHER}. We also use the notation that
$\Phi(n)$ is the number of $i$ between $1$ and $n$ inclusive, which are
relatively prime to $n$.

{\bf The Unitary Groups} Let $v_{U,n}$ be the average over
$U(n,q)$ of the order of an element of $U(n,q)$. It will be shown that for fixed
$q$ and growing $n$, $log(v_{U,n}) \geq nlog(q) - log(n) + o(log(n))$.
	
	From basic field theory the roots of an irreducible polynomial
$\phi$ of degree $n$ over $F_{q^2}$ are an orbit of some $\beta$ in a
degree $n$ extension over $F_{q^2}$ under the Frobenius map $x
\rightarrow x^{q^2}$. The next two lemmas are useful.

\begin{lemma} \label{Upolyodd} Let $L$ be a degree $n$ extension of
$F_{q^2}$, where $n$ is odd. Then an element $\beta$ of order $q^n+1$
in the multiplicative group of $L$ corresponds to an irreducible
polynomial $\phi$ of degree $n$ such that $\phi=\tilde{\phi}$.
\end{lemma}

\begin{proof}
	First note that the irreducible polynomial $\phi$ which
$\beta$ gives rise to has degree $n$. Indeed, suppose that $\beta$
lies in some proper subfield $K$ of $L$. Let $c$ be the extension
degree of $K$ over $F_{q^2}$. Then $q^n+1|q^{2c}-1$ and $c<n$. However
$c|n$, since $L$ contains $K$. This is a contradiction.

	Next, we show that $\phi=\tilde{\phi}$. Lemma $\ref{pro}$
implies that the roots of $\tilde{\phi}$ are $(\frac{1}{\beta})^
{q^{2i+1}}$ where $0 \leq i \leq n-1$. Taking $i=\frac{n-1}{2}$ shows
that $\beta$ is a root of $\tilde{\phi}$, so that $\phi=\tilde{\phi}$.
\end{proof}

\begin{lemma} \label{Upolyeven} Let $L$ be a degree $n$ extension of
$F_{q^2}$, where $n$ is even. Then an element $\beta$ of order $q^n-1$
in the multiplicative group of $L$ corresponds to an irreducible
polynomial $\phi$ of degree $\frac{n}{2}$ such that $\phi \neq
\tilde{\phi}$.  \end{lemma}

\begin{proof}
	Since the order of $\beta$ is $q^n-1$, the smallest extension
of $F_{q^2}$ containing $\beta$ is of degree $\frac{n}{2}$. Thus the
irreducible polynomial $\phi$ which $\beta$ gives rise to has degree
$\frac{n}{2}$.

	Suppose to the contrary that $\phi=\tilde{\phi}$. Then
$\beta=(\frac{1}{\beta})^{q^{2i+1}}$ for some $i$ between $0$ and
$\frac{n}{2}-1$. Thus the order of $\beta$ divides $q^c+1$ for some
$c$ between $1$ and $n-1$. This is a contradiction, so that $\phi \neq
\tilde{\phi}$.
\end{proof}

\begin{theorem} For fixed $q$ and growing $n$, $log (v_{U,n}) \geq nlog(q)
- log(n) + o(log(n))$.
\end{theorem}

\begin{proof}
	If $n$ is odd, then there are $\frac{\Phi(q^n+1)}{n}$
irreducible polynomials $\phi$ of degree $n$ such that
$\phi=\tilde{\phi}$ and the associated $\beta$ has order $q^n+1$. This
follows from Lemma $\ref{Upolyodd}$ and the fact that the elements of
order $q^n+1$ in $L$ are precisely the $\Phi(q^n+1)$ generators of the
unique order $q^n+1$ cyclic subgroup of the multiplicative group of
$L$.

	Any $\alpha$ with such a $\phi$ as its characteristic
polynomial has order $q^n+1$ and by Theorem $\ref{Uchar}$, the number
of unitary matrices with such a characteristic polynomial $\phi$ is:

\[ \frac{|U(n,q)|}{q^n+1} \]

	Thus for $n$ odd, $v_n \geq \frac{\Phi(q^n+1)}{n}$. From Stong
$\cite{St2}$ $log (\Phi(N))=log(N) + O(log(log(N)))$, so that $log (v_n)
\geq nlog(q)- log(n) + O(log(log(n)))$.

	Suppose that $n$ is even. By Lemma $\ref{Upolyeven}$ there
are $\frac{\Phi(q^n-1)} {\frac{n}{2}}$ polynomials $\phi$ of degree
$\frac{n}{2}$ such that $\phi \neq \tilde{\phi} $ and the associated
$\beta$ has order $q^n-1$. Thus there are $\frac{\Phi(q^n-1)}{n}$ such
pairs $\phi \tilde{\phi}$.  Any $\alpha$ with $\phi \tilde{\phi}$ as
its characteristic polynomial has order $q^n-1$ and by Theorem
$\ref{Uchar}$, the number of unitary matrices with such a
characteristic polynomial is:

\[ \frac{|U(n,q)|}{q^n-1} \]

	So the result follows as in the case of $n$ odd.
\end{proof}

{\bf The Symplectic Groups} Let $v_{Sp,2n}$ be the average over $Sp(2n,q)$ of the
order of an element of $Sp(2n,q)$. It will be shown that for fixed
$q$ and growing $n$, $log (v_{Sp,2n}) \geq n log(q)-log(n)+o(log(n))$. The
approach is similar to that used for the unitary groups.

	Recall from field theory that an irreducible polynomial of
degree $n$ over $F_q$ corresponds to the orbit of some $\beta$ in a
degree $n$ extension over $F_q$ under the Frobenius map $x \rightarrow
x^q$.

\begin{lemma} \label{speven} Let $L$ be a degree $2n$ extension of
$F_q$. Then an element $\beta$ of order $q^n+1$ in the multiplicative
group of $L$ corresponds to an irreducible polynomial $\phi$ of degree
2n such that $\phi=\bar{\phi}$.
\end{lemma}

\begin{proof}
	Note that the irreducible polynomial $\phi$ which $\beta$
gives rise to has degree $2n$. Suppose to the contrary that $\beta$
lied in $K$, a proper subfield of $L$.  Letting $c$ denote the
extension degree of $K$ over $F_q$, we have that $q^n+1|q^c-1$, where
$c|2n$ and $c<2n$. This is a contradiction.

	By Lemma $\ref{product2}$, the roots of $\bar{\phi}$ are
$(\frac{1}{\beta})^{q^i}$ where $1 \leq i \leq 2n$. Taking $i=n$ shows
that $\beta$ is a root of $\bar{\phi}$. Thus $\phi=\bar{\phi}$.
\end{proof}

	Let $\Phi(n)$ be the number of $i$ between 1 and $n$ inclusive
which are relatively prime to $n$.

\begin{theorem} \label{ordsym} For fixed $q$ and growing $n$,
$log(v_{Sp,2n}) \geq nlog(q) - log(n)+ o(log(n))$.  \end{theorem}
	
\begin{proof}
	From Lemma $\ref{speven}$, there are $\frac{\Phi(q^n+1)}{2n}$
irreducible polynomials $\phi$ of degree $2n$ such that the associated
$\beta$ has order $q^n+1$. This is because the elements of order $q^n+1$ in
$L$ are the $\Phi(q^n+1)$ generators of the order $q^n+1$ cyclic subgroup
of the multiplicative group of $L$.

	By Theorem $\ref{SpChar}$, the number of elements of $Sp(2n,q)$
with characteristic polynomial $\phi$ is:

\[ \frac{|Sp(2n,q)|}{q^n+1} \]

	So considering only such $\alpha$ gives the lower bound
$v_{Sp,2n} \geq \frac{\Phi(q^n+1)}{2n}$. Using the fact from Stong
$\cite{St2}$ that $log(\Phi(N))= log(N)+O(log(log(N)))$ proves the
theorem.
\end{proof}

{\bf The Orthogonal Groups}

	Let $v_{O,n}$ be one half of the sum of the average orders of
elements of $O^+(n,q)$ and $O^-(n,q)$.

\begin{theorem} For fixed $q$ and growing $n$, $log(v_{O,n}) \geq
\frac{n}{2}log(q) - log(n) + o(log(n))$.
\end{theorem}

\begin{proof}
	Assume that $n$ is even (the case of $n$ odd being similar),
and consider only orthogonal matrices $\alpha$ whose characteristic
polynomial comes from a $\beta$ of order $q^{\frac{n}{2}}+1$ in a
degree $n$ extension of $F_q$. Then use Lemma $\ref{speven}$ and
Theorem $\ref{Ocharpoly}$ and argue as in Theorem $\ref{ordsym}$.
\end{proof}

\section{Suggestions for Future Research} \label{SUGGESTIONS}

\begin{enumerate}

\item Read information off of the cycle indices for the classical
groups. For instance extend Goh and Schmutz's $\cite{Go2}$ theorem that the
number of Jordan blocks is asymptotically normal to the other classical
groups. Give upper bounds for the average order of an element of a unitary,
symplectic, or orthogonal matrix. Carry over the asymptotic work on
semisimple, regular, and regular-semisimple elements in the general linear groups.

\item Study the cycle indices for the characteristic 2 cases for the
symplectic and orthogonal groups. One can deduce from Wall
$\cite{Wal}$ that these also factor and that results analogous to
Theorem $\ref{Count2}$ go through.

\item Find Lie Algebra cycle indices. For instance Stong's cycle
index for $Mat(n,q)$ (see Section $\ref{GLCYC}$) encodes the
orbits of the adjoint action of $GL(n,q)$ on its Lie Algebra. Are
there factorizations for the other classical groups as well?

\end{enumerate}

\section{Acknowledgements}

		This work is taken from the author's Ph.D. thesis, done under the
supervision of Persi Diaconis at Harvard University. The author is deeply indebted
to him for emphasizing the importance of cycle indices and much helpful advice.
The author also thanks Dick Gross for patiently explaining the basics of algebraic
groups. This research was done
under the generous 3-year support of the National Defense Science and Engineering
Graduate Fellowship (grant no. DAAH04-93-G-0270) and the support of the Alfred P.
Sloan Foundation Dissertation Fellowship.

\end{document}